\documentclass[review]{elsarticle}
\usepackage{color}
\usepackage{lineno,hyperref,amsfonts,amsmath}
\modulolinenumbers[5]

\journal{Journal of \LaTeX\ Templates}

\usepackage{cancel}
\newcommand{\norm}[1]{\left\Vert#1\right\Vert}
\newtheorem{lemma}{Lemma}[section]
\newtheorem{remark}{Remark}[section]

\newtheorem{theorem}{Theorem}[section]

\newtheorem{definition}{Definition}[section]

\newtheorem{proposition}{Proposition}[section]
\DeclareMathOperator{\sgn}{sgn}

\def\={\buildrel \triangle \over =}

\begin{document}

\begin{frontmatter}

\title{Non-homogeneous  initial boundary value problems for the biharmonic  Schr\"odinger equation on an interval}


\author[mysecondaryaddress]{Junfeng Li}
\ead{ junfengli@dlut.edu.cn}
\address[mysecondaryaddress]{School of Mathematical Sciences, Dalian University of Technology, Dalian, L.N. China.  }

\author[mysecondaryaddress]{Chuang Zheng\corref{mycorrespondingauthor}}
\cortext[mycorrespondingauthor]{Corresponding author}
\ead{chuang.zheng@bnu.edu.cn}

\address[mymainaddress]{School of Mathematical Sciences,
	Beijing Normal University, Beijing,	China. }

\begin{abstract}
In this paper we consider the initial boundary value problem (IBVP) for the nonlinear biharmonic Schr\"odinger equation posed on a bounded interval $(0,L)$ with non-homogeneous  Navier or Dirichlet  boundary conditions, respectively. For Navier boundary IBVP, we set up its local well-posedness if the initial data lies in $H^s(0, L)$  with $s\geq 0$ and $s\neq n+1/2, n\in \mathbb{N}$, and the boundary data are selected from the appropriate spaces with optimal regularities, i.e., the $j$-th order data are chosen in  $H_{loc}^{(s+3-j)/4}(\mathbb {R}^+)$, for $j=0,2$.  For Dirichlet boundary IBVP the corresponding local well-posedness is obtained when $s>10/7$ and $s\neq n+1/2, n\in \mathbb{N}$, and the boundary data are selected from the appropriate spaces with optimal regularities, i.e., the $j$-th order data are chosen in  $H_{loc}^{(s+3-j)/4}(\mathbb {R}^+)$, for $j=0,1$. 
\end{abstract}


\end{frontmatter}


\section{Introduction}

Biharmonic Schr\"odinger equations have been introduced in many scientific fields such as quantum mechanics, nonlinear optics and plasma physics. For instance, its nonlinear form is given in \cite{MR1781877} to study  the motion of a vortex filament in an incompressible fluid  and also for the small fourth order dispersion terms in the propagation of intense laser beams in a bulk medium with Kerr nonlinearity (\cite{Karpman96} and \cite{MR1779828}).   The well-posedness and other properties of the solutions have been intensively studied from the mathematical point of view in the whole domain $\mathbb{R}^n$  and $\mathbb{T}^n$,  respectively(see \cite{ MR1745182, MR3665549, MR2230468, MR3542971, MR2353631, MR2502523,  MR1696311}  and the references therein). As far as we can see in most of the existing references, the boundary data are neglected  due to the fact that the system under consideration evolves in the unbounded domain. In some cases,  the system will be required to evolve in a finite domain,  such as implementing the numerical simulation, or imposing controls on the boundary, etc. See, for instance, \cite{MR3720366, MR3356487, MR2928966}. In these cases,  it is reasonable to take the boundary data not only with the homogeneous Dirichlet/Neumann boundary conditions but also with  non-homogeneous boundary conditions.  The monographs \cite{MR0350177, MR0350178, MR0350179} provide a systematical way to analyze the abstract model of the non-homogeneous boundary value problems via PDE techniques. Recently,  in \cite{MR1998942}, by means of the boundary integral operator and harmonic analysis,  the well-posedness of the  nonlinear KdV equation on a bounded interval  is established and the regularities of the boundary data are given.  Consequently, several works have been done for the second order Schr\"odinger equation, KdV and Kuramoto-Sivashinsky equation, etc. (\cite{MR3734975,  CSZ, MR3702485, MR3913121, MR3809532}).

Compared to the second order system,  the biharmonic operator systems  serve a rich varieties of choices on the boundary data and corresponds to different physical phenomena when the systems evolve in a finite domain. For instance, in the case of biharmonic plate model,  the Dirichlet conditions correspond to the clamped plate model and  the Navier or Steklov boundary data correspond to the hinged plate model, either by neglecting or considering the contribution of the curvature of the boundary. Each of these boundary data has its own feature and leads to different physical properties.  The aim of the present work is to analyze the initial-boundary-value-problems (IBVP henceforth) of the biharmonic Schr\"odinger equations with the typical Navier and Dirichlet boundary conditions.

More precisely,  we consider the IBVP of the following biharmonic nonlinear Schr\"odinger equation posed on the finite interval $(0,L)$, i.e.
\begin{equation}\label{IBVP}
i\partial_t u+ \partial_x^4u+\lambda |u|^{p-2}u=0, \qquad x\in (0,L), \; t\in {\mathbb {R}}
\end{equation}
with initial data
\begin{equation}\label{ID}
u(x,0)=\varphi(x), \qquad x\in (0,L).
\end{equation}
Here the parameter $\lambda$ is a non-zero real number and $p\geq 3$.  The  Navier boundary data  are described by (\cite{MR2667016})
\begin{equation}\label{NBC}
u(0,t)=h_1(t), u(L,t)=h_2(t), \partial_x^2 u(0,t)= h_5(t), \partial_x^2 u(L,t)=h_6(t),\;\; t\in {\mathbb {R}}.
\end{equation}
The  Dirichlet boundary data  are described by
\begin{equation}\label{DBC}
u(0,t)=h_1(t), u(L,t)=h_2(t), \partial_x u(0,t)= h_3(t), \partial_x u(L,t)=h_4(t),\;\; t\in {\mathbb {R}}.
\end{equation}

We discuss the well-posedness of equation \eqref{IBVP} corresponding to specific initial data in the $L^2$-based Sobolev spaces $H^s(0,L)$ with  boundary data \eqref{NBC} or \eqref{DBC} . We assume that the compatibility conditions
$$
h_1(0)=\varphi(0), h_2(0)=\varphi(L),  h_3(0)=\varphi_{x}(0), h_4(0)=\varphi_{x}(L),
$$
are valid if $s>\frac{3}{2}$ for \eqref{DBC}.  Moreover,
$$
h_1(0)=\varphi(0), h_2(0)=\varphi(L),  h_5(0)=\varphi_{xx}(0), h_6(0)=\varphi_{xx}(L),
$$
are valid if $s>\frac{5}{2}$ for \eqref{NBC}. If $s>0$ is large, we also assume that $|u|^{p-2}u$ is differentiable. Furthermore, for the convenience of our discussion on the trace of functions in $H^s(\mathbb{R})$, we always implicitly assume that
$$
s\neq n+\frac{1}{2}, \qquad \hbox{for} \qquad n=0,1,2,\cdots.
$$

We extend the approach for studying IBVP of the second order Schr\"odinger equation and the Kuramoto-Sivashinsky equation (\cite{MR3734975,MR3702485}) and obtain two local well-posedness theorems. More precisely,  we first establish the result with Navier boundary conditions:
\begin{theorem}\label{mainN}
	Assume that $s>\frac{1}{2}$ with $p\in[3, \infty)$ and $T>0$. Let
	$\varphi\in H^s(0,L)$, $(h_1,  h_5)$ and $(h_2, h_6)$ be in the space $X_N$ with
	$$
	X_N:=  H^{\frac{s+3}{4}}(0,T)\times H^{\frac{s+1}{4}}(0,T).
	$$
	The IBVP \eqref{IBVP}--\eqref{NBC} is locally well-posed in $H^s(0,L)$.
	
	If $s\in [0, \frac{1}{2})$ with $p\in[3, 4]$, IBVP \eqref{IBVP}--\eqref{NBC} is  locally well-posed in $H^s(0,L)$  for $(\varphi, h_1, h_2,  h_5, h_6)$ being in the same spaces.
	
	More precisely,  for any $\gamma>0$, there exists $T^*$ with $T^*\in(0,T]$ depending only on $s, \gamma$ and $T$ such that if
	$$
	\norm{\varphi}_{H^s(0,L)}+\norm{(h_1,  h_5)}_{X_N}+\norm{(h_2, h_6)}_{X_N}\leq \gamma
	$$
	the IBVP \eqref{IBVP}--\eqref{NBC} admits a unique solution $u\in C([0, T^*]; H^s(0,L))$.
\end{theorem}

The second result is about the result  with the Dirichlet boundary conditions:
\begin{theorem}\label{mainD}
	Assume that $s>\frac{10}{7}$ with $p\in[3, \infty)$ and $T>0$. Let
	$\varphi\in H^s(0,L)$, $(h_1,  h_3)$ and $(h_2, h_4)$ be in the space $X_D$ with
	$$
	X_D:=  H^{\frac{s+3}{4}}(0,T)\times H^{\frac{s+2}{4}}(0,T).
	$$
	The IBVP \eqref{IBVP}, \eqref{ID}, \eqref{DBC} is locally well-posed in $H^s(0,L)$.
	
	More precisely,  for any $\gamma>0$, there exists  $T^*$ with $T^*\in(0,T]$ depending only on $s, \gamma$ and $T$ such that if
	$$
	\norm{\varphi}_{H^s(0,L)}+\norm{(h_1,  h_3)}_{X_D}+\norm{(h_2, h_4)}_{X_D}\leq \gamma
	$$
	the IBVP \eqref{IBVP}, \eqref{ID}, \eqref{DBC}  admits a unique solution $u\in C([0, T^*]; H^s(0,L))$.
\end{theorem}

Some remarks are in order:

\begin{remark}\label{Kato-Smooth}
	Note that the regularity of the boundary data is related to the Kato smoothing property which is extensively studied in the case with unbounded domain (see, for instance, \cite{MR3577868}). In fact, the stipulation $s'=\frac{1}{2}(s+\frac{1}{2})$ is fulfilled for  the second order problem posed on the half line  $\mathbb{R}^+$. Here, $s$  and $s'$ represent the regularity of the initial data  and the optimal regularity of the zero-order boundary data, respectively.  This is the exact ``cut-off " regularity predicted by Kato smoothing property.  However, in the case with the bounded domain $(0,L)$, the optimal relation  jumps to $s'=\frac{1}{2}(s+1)$, which is surprisingly higher than one expected to be (\cite{MR3734975}).
	Our results in this paper offer a clear explanation of this  difference: The Kato smoothing property has been {\bf{doubled}} while the dispersive equations evolve on the bounded domain, i.e.  the following formula holds:
	\begin{equation}\label{KSB}
	s_i'=\frac{1}{4}\Big (s+2\frac{4-1}{2}-i\Big ) \qquad \hbox{for} \qquad i=0, 1, 2.
	\end{equation}
	In \eqref{KSB}, $s_{i}'$ is the optimal regularity of the $i$-th order boundary data for $i=0,1, 2$, respectively. See Lemma \ref{opt0} for more details. 
\end{remark}

\begin{remark}
	Note that system \eqref{IBVP}, \eqref{ID} with Navier boundary conditions \eqref{NBC}  is in accordance to the case of Dirichlet boundary conditions for the second order Schr\"odinger equation.  
	In both cases, the solutions evolved on $(0, L)$ can be seen as the restriction of  solutions of IVP (p35 of \cite{MR3734975} and \eqref{wn1} in our paper) which evolves on $\mathbb{R}$ with odd initial data generating from $(0, L)$.  Note that the result can also be extended to the case with boundary conditions
	$$
	\partial_x u(0,t)=h_3(t), \partial_x u(L,t)=h_4(t), \partial_x^3 u(0,t)= h_7(t), \partial_x^3 u(L,t)=h_8(t),\;\; t\in {\mathbb {R}}
	$$
	and the optimal regularities will follow the formula \eqref{KSB} with $i=1$ and $i=3$.
\end{remark}
\begin{remark}
	Note that the same methodology for the case with  Navier boundary conditions does not fit for the case with Dirichlet boundary conditions of the biharmonic operator.  In fact, the solution of \eqref{IBVP}, \eqref{ID} with \eqref{DBC}   contains two kinds of solutions of IVP, which are evolving on $\mathbb{R}$ with odd and even initial data generating from $(0, L)$, respectively. We offer a new idea by introducing a special Fourier expansion \eqref{phioe} and, consequently,  Dirichlet boundary data can be transformed to boundary integrals, as the form in Proposition \eqref{BoIn4}. The details of the proof are put in Section \ref{3.1}.  
\end{remark}

\begin{remark}For IBVPs, it is crucial to consider the effects on boundary integrals raising by the initial data and boundary data, respectively. These effects for the Navier problem can be handled by odd extension of the equation from the interval $[0,1]$ into $[-1,1]$. However, it does not work for the case  with Dirichlet boundary data. For the latter case, we need to carefully check whether the boundary values raised by the initial data can share the same regularity of those raised by the non-homogeneous  boundary data (see, for instance, the extra terms $r_i, i=1,2,3,4$ in Lemma \ref{lema1} and Lemma \ref{lema2}). Consequently, the regularity of the initial data can not be too low. This effect brings us the  lower bound $s>10/7$ in Theorem \ref{mainD}. It is still an open problem whether it is sharp. 
	
\end{remark}

The rest of the paper is organized as follows.  In Section \ref{prelim}, we  collect some basic definitions and give a rule between the optimal regularities of the boundary data and the regularity of the initial data.  Section \ref{3.1} has two subsections. We establish the estimates of the boundary integrals which coincide to the Navier and Dirichlet non-homogeneous boundary data, respectively.  The  well-posedness of the corresponding nonlinear problems are given as two subsections accordingly in Section \ref{3.2}.  We put some technical details  in  Appendix and construct  counterexamples to verify the optimality of the regularities of the boundary data.

\section{Preliminary}\label{prelim}


We first state a precise definition of well-posedness for the problem \eqref{IBVP}, \eqref{ID} with Navier boundary data \eqref{NBC}.
\begin{definition}\label{Navier}
	The IBVP \eqref{IBVP} -- \eqref{NBC} is said to be (locally) well-posed in $H^s(0,L)\times H^{s_0'}(0,T)\times H^{s_2'}(0,T)$  if for $\varphi\in H^s(0,L)$, $h_1, h_2\in H^{s_0'}(0,T)$ and  $ h_5, h_6\in H^{s_2'}(0,T)$ satisfying certain natural compatibility conditions, there exists a $T'\in(0,T]$ depending only on 
	$$
	r=\norm{\varphi}_{H^{s}(0,L)}+\sum_{j=1,2}\Big (\norm{h_{j}}_{H^{s_0'}(0,T)}+\norm{h_{j+4}}_{H^{s_2'}(0,T)}\Big )
	$$
	such that \eqref{IBVP} -- \eqref{NBC} admits a unique solution $u\in  C([0,T];H^s(0,L))$.  Moreover, the solution depends continuously on $(\varphi, h_1,h_2, h_5,h_6)$ in the corresponding spaces.
\end{definition}
The corresponding solutions  of the system are  defined by
\begin{definition}
	Let $s\leq 4$, $s_2'\leq s_0'\leq s$ and $T>0$ be given. For any $\varphi\in H^s(0,L)$, $h_1, h_2\in H^{s_0'}(0,T)$ and $ h_5, h_6\in H^{s_2'}(0,T)$, we say that $u\in C([0,T];H^s(0,L))$ is a solution of \eqref{IBVP} -- \eqref{NBC},  if there exists a sequence
	$$
	u_n\in C([0,T];H^4(0,L))\cap C^1([0,T];L^2(0,L)), \quad n=1,2,3,\cdots
	$$
	such that
	\begin{enumerate}[1)]
		\item $u_n$ satisfies the equation of \eqref{IBVP} in $L^2(0,L)$ for $0\leq t\leq T$,
		\item $u_n$ converges to $u$ in $C([0,T];H^s(0,L))$ as $n\rightarrow\infty$,
		\item $\varphi_n(x)=u_n(x,0)$ converges to $\varphi(x)$ in $H^s(0,L)$ as $n\rightarrow\infty$,
		\item $h_{1,n}(t)=u_n(0,t), h_{2,n}(t)=u_n(L,t)$ are in $H^{s_0'}(0,T)$ and converge to $h_1(t)$ and $h_2(t)$, respectively, in $H^{s_0'}(0,T)$ as $n\rightarrow\infty$,
		\item $h_{5,n}(t)= \partial_x^2 u_{n}(0,t),  h_{6,n}(t)=  \partial_x^2 u_{n}(L,t)$ are in $H^{s_2'}(0,T)$ and converge to $ h_5(t)$ and $h_6(t)$, respectively, in $H^{s_2'}(0,T)$ as $n\rightarrow\infty$.
	\end{enumerate}
\end{definition}

Similarly, we give a precise definition of well-posedness for the problem \eqref{IBVP}, \eqref{ID} with Dirichlet boundary data \eqref{DBC}.
\begin{definition}\label{Dirichlet}
	The IBVP \eqref{IBVP}, \eqref{ID}, \eqref{DBC} is said to be (locally) well-posed in $H^s(0,L)\times H^{s_0'}(0,T)\times H^{s_1'}(0,T)$  if for $\varphi\in H^s(0,L)$, $h_1, h_2\in H^{s_0'}(0,T)$ and  $ h_3, h_4\in H^{s_1'}(0,T)$ satisfying certain natural compatibility conditions, there exists a $T'\in(0,T]$ depending only on 
	$$
	r=\norm{\varphi}_{H^{s}(0,L)}+\sum_{j=1,2}\Big (\norm{h_{j}}_{H^{s_0'}(0,T)}+\norm{h_{j+2}}_{H^{s_1'}(0,T)}\Big )
	$$
	such that \eqref{IBVP}, \eqref{ID},  \eqref{DBC} admits a unique solution $u\in  C([0,T]; H^s(0,L))$.  Moreover, the solution depends continuously on $(\varphi, h_1,h_2, h_3,h_4)$ in the corresponding spaces.
\end{definition}
The corresponding solutions  of the system are  defined by
\begin{definition}
	Let $s\leq 4$, $s_1'\leq s_0'\leq s$ and $T>0$ be given. For any $\varphi\in H^s(0,L)$, $h_1, h_2\in H^{s_0'}(0,T)$ and $ h_3, h_4\in H^{s_1'}(0,T)$, we say that $u\in C([0,T];H^s(0,L))$ is a solution of \eqref{IBVP}, \eqref{ID} and  \eqref{DBC},  if there exists a sequence
	$$
	u_n\in C([0,T];H^4(0,L))\cap C^1([0,T];L^2(0,L)), \quad n=1,2,3,\cdots
	$$
	such that
	\begin{enumerate}[1)]
		\item $u_n$ satisfies the equation of \eqref{IBVP} in $L^2(0,L)$ for $0\leq t\leq T$,
		\item $u_n$ converges to $u$ in $C([0,T];H^s(0,L))$ as $n\rightarrow\infty$,
		\item $\varphi_n(x)=u_n(x,0)$ converges to $\varphi(x)$ in $H^s(0,L)$ as $n\rightarrow\infty$,
		\item $h_{1,n}(t)=u_n(0,t), h_{2,n}(t)=u_n(L,t)$ are in $H^{s_0'}(0,T)$ and converge to $h_1(t)$ and $h_2(t)$, respectively, in $H^{s_0'}(0,T)$ as $n\rightarrow\infty$,
		\item $h_{3,n}(t)= \partial_x u_{n}(0,t),  h_{4,n}(t)=  \partial_x u_{n}(L,t)$ are in $H^{s_1'}(0,T)$ and converge to $ h_3(t)$ and $h_4(t)$, respectively, in $H^{s_1'}(0,T)$ as $n\rightarrow\infty$.
	\end{enumerate}
\end{definition}



Now we discuss the relationship between $s_0', s_1', s_2'$ and $s$ in the definition of well-posedness. Note that the subscript $i$ of $s_i'$ corresponds to the $i$-th order spatial derivative on the boundary. For the linear Schr\"odinger equation on the whole line
$$
i\partial_t   v+\partial_x^4 v=0, \qquad v(x,0)=\varphi(x), \qquad\hbox{for}\quad x\in\mathbb{R},
$$
the Kato smoothing property affirms that $\varphi\in H^s(\mathbb{R})$ implies $v\in L^2_{loc}(\mathbb{R}, H_{loc}^{s+\frac{3}{2}}(\mathbb{R}))$
\footnote{In \cite{MR3577868}, it is shown that for the general dispersive equations
	$$
	\omega_t+i P(D)\omega=0, \quad \omega(x,0)=q(x), \quad x\in \mathbb{R}, \quad t\in \mathbb{R}
	$$	
	enjoy the sharp Kato smoothing property
	$$
	q\in H^s(\mathbb{R})\Longrightarrow \omega\in L^2_{loc}\Big (\mathbb{R}; H^{s+\frac{n-1}{2}}_{loc}(\mathbb{R})\Big )
	$$
	where $n$ is the order of the pseudo-differential operator $P(D)$.
}(\cite{MR928265, MR3577868}).
On the other hand, the biharmonic operator suggests the relation $\partial_t \sim \partial_{x}^4$. Combining these two facts one can easily find that, for the half line domain, the formula for boundary regularities must be
$$
s_0'=\frac{1}{4}\Big (s+\frac{3}{2}\Big ).
$$
This methodology works well for the second order Schr\"odinger equation (see  (2.12) of \cite{MR3734975}).

However, in the case of  the bounded interval $(0,L)$,  the correct value of the sharp regularity is higher (see, for instance, (2.13) of \cite{MR3734975}).
Based on the result in this paper, one could expect that for the $2m$-th order Schr\"odinger equation, they obey the  stipulation
$$
s_0'=\frac{1}{2m}\Big (s+2\frac{2m-1}{2}\Big ),
$$
where $\frac{2m-1}{2}$ is the effect of the Kato smoothing.

In addition, if we denote by $s_{i}'$ the optimal regularity of the $i$-th order boundary data, i.e.,
$$
\partial_x^i u(0, t)\in H_{loc}^{s_i'}(\mathbb{R}), \qquad \hbox{for} \qquad i=0,1,\dots, 2m-1,
$$
the relation between $s$ and $s_i'$ should be
$$
s_i'=\frac{1}{2m}\Big (s+2\frac{2m-1}{2}-i\Big ) \qquad \hbox{for} \qquad i=0,1,\dots, 2m-1.
$$
So far, the above rule is correct for $m=1, i=0$(\cite{MR3734975}),  $m=2, i=0, 2$ (with Navier boundary data) and  $m=2, i=0, 1$ (with Dirichlet boundary data). As an ongoing problem, it would be stirring to establish a uniform proof for arbitrary integer $m$.

\section{The linear problem}\label{3.1}

We analyze the linear fourth order Schr\"odinger equation with two types of boundary conditions, i.e. with Navier boundary conditions and with Dirichlet boundary conditions in Subsection \ref{3.11} and Subsection \ref{3.12}, respectively.

Without loss of generality, we take $L=1$.  

\subsection{Case 1: With Navier boundary conditions}\label{3.11}

To begin with, consider the IBVP
\begin{equation}\label{LinearIBVP}
\left\{\begin{array}{ll}
i\partial_t u+  \partial_x^4u=0, & x\in (0,1), \; t\in {\mathbb {R}}\\
u(x,0)=\varphi(x),&x\in (0,1)\\
u(0,t)=u(1,t)=0,\quad &t\in{\mathbb {R}} \\
\partial_x^2 u(0,t)=  \partial_x^2 u(1,t)=0, &t\in{\mathbb {R}}
\end{array}\right.
\end{equation}
for the linear Schr\"odinger equation. According to the standard semi-group theory, for any $\varphi\in L^2(0,1)$, the IBVP admits a unique solution $u\in C({\mathbb {R}}^+; L^2(0,1))$ given by
\begin{equation}\label{wnphio}
u(t)=W^N(t)\varphi
\end{equation}
where $W^N(t)$ is the $C_0$-group in $ L^2(0,1)$ generated by the operator $A^Nv=i v''''$ with domain $\mathcal{D}(A^N)=\{u\in H^4(0,1): u(x)=u_{xx}(x)=0, \; x=0, 1\} $.

In fact, due to the homogeneous boundary data, the solution $u$ can be expressed by Fourier sine series, i.e.,
\begin{equation}\label{linear}
u(x,t)=W^N(t)\varphi:=\sum_{n=1}^{\infty} \hat{\varphi}(n)e^{i(n\pi)^4t}\sin{n\pi x},
\qquad  x\in(0,1),
\end{equation}
where $\hat{\varphi}(n)=2\int_0^1\varphi(x)\sin{n\pi x} dx$.  This can be written in the complex form
$$
u(x,t)=\sum_{n\in \mathbf{Z}} e^{i(n\pi)^4t}\hat{\phi}(n)e^{in\pi x}
$$
where
$$
\hat{\phi}(n)=
\left\{\begin{array}{ccc}
&\hat{\varphi}(n)   & n\geq1,\\
&0  & n=0,\\
&-\hat{\varphi}(n)  & n\leq-1.
\end{array} \right.
$$
In this form, $u$ may be viewed as  the solution of the corresponding Cauchy problem with periodic domain, i.e.
\begin{equation}\label{perio-2}
\left\{\begin{array}{ll}
i\partial_t u+  \partial_x^4u=0, & x\in (-1,1), \; t\in {\mathbb {R}}\\
u(x,0)=\varphi^*,& x\in (-1,1)\\
u(-1,t)=u(1,t),\;
\partial_x u(-1,t)=  \partial_x u(1,t), &t\in{\mathbb {R}}\\
\partial_x^2 u(-1,t)=\partial_x^2 u(1,t),\;
\partial_x^3 u(-1,t)=  \partial_x^3 u(1,t), &t\in{\mathbb {R}}
\end{array}\right.
\end{equation}
If $u$ is a solution of \eqref{perio-2} with odd initial data $\varphi^*$, it is obvious that its restriction to $(0,1)$ solves the corresponding linear problem \eqref{LinearIBVP}. Thus
\begin{equation}\label{wn1}
[W_{\mathbb{T}}\varphi^*](x)=[W^N\varphi](x), \qquad x\in(0,1),
\end{equation}
where $W_{\mathbb{T}}$ is the $C_0$-group in $L^2(\mathbb{T})$ generated by the operator $A_{\mathbb{T}}$ in $L^2(\mathbb{T})$.
We first give the $L^4$ estimate of \eqref{LinearIBVP} by Bourgain's theory:
\begin{proposition}\label{L4estimates}
	Let $s\geq 0$ and $T>0$  be given and let $\Omega_T=(0,1)\times(0,T)$. For any $\varphi\in H^s(0,1), u=W^N(t)\varphi\in L^4(\Omega_T)\cap C([0,T]; H^s(0,1))$ has
	$$
	\norm{u}_{L^4(\Omega_T)\cap C([0,T]; H^s(0,1))}\leq C\norm{\varphi}_{H^s(0,1)},
	$$
	where $C$ depends only on $s$ and $T$.
\end{proposition}
\textbf{ Proof:} No loss the generalization, we may assume the time existence to be $\mathbb T=[-1,1]$. Then we  denote $\Omega_T=\mathbb T^2$.
We  only need to prove
$$
\norm{W_{\mathbb T}\varphi^*}_{L^4(\mathbb T^2)}\leq C\norm{\varphi^*}_{L^2(\mathbb T)}.
$$
For $s\geq0$, the estimates follow from Plancherel theorem and the embedding property of $H^s(0,1)$. We write
$$
\norm{\sum_{k\in \mathbb{Z}}\hat{\phi}(k) e^{i(\pi kx+(\pi k)^4t)}}_{L^4({\mathbb {T}}^2)}=\norm{\sum_{k,l\in\mathbb{Z}}\hat{\phi}(k)\bar{\hat{\phi}}(l) e^{i[\pi(k-l)x+(\pi^4(k^4-l^4))t]}}_{L^2({\mathbb {T}}^2)}.
$$
By Plancherel theorem, it equals to
$$\Bigg\{\sum_{\xi,\eta\in\mathbb{Z}}\Big(\sum_{(k,l)\in A(\xi,\eta)}\hat{\phi}(k)\bar{\hat{\phi}}(l)\Big)^2\Bigg\}^{\frac{1}{2}}.$$
Here $A(\xi,\eta)=\{(k,l)\in\mathbb{Z}^2;k-l=\xi,k^4-l^4=\eta\}.$ It is not hard to prove that $A(\xi,\eta)\cap A(\xi^\prime,\eta^\prime)=\emptyset.$ Meanwhile, for any fixed $\xi,\eta\in\mathbb{Z}$,
$$\#A(\xi,\eta)\leq 3.$$
We finish the proof by Cauchy-Schwarz.

Next, we consider the corresponding non-homogeneous problem
$$
\left\{\begin{array}{ll}
i\partial_t u+  \partial_x^4u=f, & x\in (0,1), \; t\in {\mathbb {R}}\\
u(x,0)=0,&x\in (0,1)\\
u(0,t)=u(1,t)=0,\quad &t\in{\mathbb {R}} \\
\partial_x^2 u(0,t)=  \partial_x^2 u(1,t)=0, &t\in{\mathbb {R}}.
\end{array}\right.
$$
By Duhamel's principle, its solution can be written as
$$
u(t)=-i\int_0^tW^N(t-\tau)f(\cdot, \tau)d\tau.
$$
\begin{proposition}\label{p4.6-0}
	Let $s\in[0,4]$ and $T>0$ be given. Let
	$$
	u(t)=W^N(t)\varphi, \qquad v(t)=\int_0^t W^N(t-\tau)f(\cdot, \tau)d\tau
	$$
	and
	$$
	w(t)=\int_0^t W^N(t-\tau)g(\cdot, \tau)d\tau
	$$
	with $\varphi\in H^s(0,1)$, $f\in L^1(0,T; H^s(0,1))$ and $g\in W^{\frac{s}{4},1}(0,T; L^2(0,1))$ satisfying
	$$
	\varphi(0)=\varphi(1)=0,\;\varphi_{xx}(0)=\varphi_{xx}(1)=0,\; f(0,t)=f(1, t)=\partial_x^2 f(0,t)=\partial_x^2 f(1,t)\equiv0
	$$
	when $s>\frac{5}{2}$ and
	$$
	\varphi(0)=\varphi(1)=0, \; f(0,t)=f(1, t)\equiv0
	$$
	with $s\in (\frac{1}{2},\frac{5}{2})$. Then $u, v, w\in C([0,T]; H^s(0,1))$ and
	\begin{eqnarray}
	\label{TPex1}\norm{u}_{C([0,T]; H^s(0,1))}&\leq & C_{T,s}\norm{\varphi}_{H^s(0,1)}\\
	\label{TPex2}\norm{v}_{C([0,T]; H^s(0,1))}&\leq & C_{T,s}\norm{f}_{L^1(0,T; H^s(0,1))}\\
	\label{TPEx}\norm{w}_{C([0,T]; H^s(0,1))}&\leq & C_{T,s}\norm{g}_{W^{\frac{s}{4},1}(0,T; L^2(0,1))}
	\end{eqnarray}
	where the constant $C_{T,s}$ depends only on $s$ and $T$.
\end{proposition}

\textbf{Proof:}  Recalling \eqref{linear}, the Plancherel theorem gives
\begin{equation}\label{Plachere}
\norm{u(t)}_{L^2(0,1)}=\norm{\varphi}_{L^2(0,1)}\qquad \hbox{for all }\quad t\in{\mathbb {R}}.
\end{equation}
Thus the first two estimates \eqref{TPex1} and \eqref{TPex2} hold for $s\in\mathbb R$. We now give the proof of \eqref{TPEx}. The case $s=0$ follows from \eqref{Plachere}. We give the proof the case $s=4$. And the case of $0<s<4$ follows from interpolation argument.
By the definition
\begin{eqnarray*}
	\norm{w(t)}_{H^4(0,1)}&=&\left(\sum_{n\in\mathbb Z} |n\pi|^8\left|\int_0^t\hat{g}(n,\tau)e^{-i(n\pi)^4(t-\tau)}d\tau\right|^2\right)^{\frac12}\\
	&=&\left(\sum_{n\in\mathbb Z} \left|\int_0^{t}\hat{g}(n,\tau)de^{i(n\pi)^4\tau}\right|^2\right)^{\frac12}\\
	&=&\left(\sum_{n\in\mathbb Z} \left|\hat{g}(n,t)e^{i(n\pi)^4t}-\hat{g}(n,0)-\int_0^{t}\hat{g}_\tau(n,\tau)e^{i(n\pi)^4\tau}d\tau\right|^2\right)^{\frac12}\\
\end{eqnarray*}
Then \eqref{TPEx} follows from Sobolev embedding and Minkowski's inequality.

Now we consider the linear problem with non-homogeneous Navier boundary conditions, i.e.
\begin{equation}\label{4th-BD-02}
\left\{\begin{array}{ll}
i\partial_t u+  \partial_x^4u=0, & (x, t)\in (0,1)\times \mathbb{R}\\
u(x,0)=0,&x\in (0,1)\\
u(0,t)=h_1(t),\quad u(1,t)=h_2(t),\quad &t\in\mathbb{R} \\
\partial_x^2 u(0,t)= h_5(t), \quad  \partial_x^2 u(1,t)=h_6(t), &t\in\mathbb{R}.
\end{array}\right.
\end{equation}
We put the compatibility conditions
$$
h_1(0)=h_2(0)=0, \qquad  h_5(0)=h_6(0)=0
$$
if it is necessary.
\begin{proposition}\label{expruxt-02} 
	The solution of \eqref{4th-BD-02} can be expressed as
	\begin{eqnarray}
	\nonumber u(x,t) &=& \sum_{k=1}^\infty 2i(k\pi)^3\int_0^t e^{i(k\pi)^4(t-\tau)}(h_1(\tau)-\cos(k\pi) h_2(\tau))d\tau \sin(k\pi x)\\
	\nonumber        & & +\sum_{k=1}^\infty (-2ik\pi)\int_0^t e^{i(k\pi)^4(t-\tau)}( h_5(\tau)-\cos(k\pi) h_6(\tau))d\tau \sin(k\pi x) \\
	\nonumber
	&\=& W_{0,N} h_1 +\left(W_{0,N} h_2\right)\left|_{x\rightarrow 1-x}\right.+W_{2,N}  h_5+\left(W_{2,N} h_6\right)\left|_{x\rightarrow 1-x}.\right.
	\end{eqnarray}
\end{proposition}
\textbf{Proof:} We consider the special case where $h_2\equiv h_6\equiv 0$ and $h_1(0)= h_5(0)=0$. We define $v$ by
\begin{equation}\label{u-vrelation}
u(x,t)=v(x,t)+(1-x)(h_1-\frac{1}{6} h_5)(t)+\frac{1}{6}(1-x)^3 h_5(t).
\end{equation}
Then $v(x,t)$ solves
$$
\left\{\begin{array}{ll}
i\partial_t v+  \partial_x^4 v=f(x,t)
& x\in (0,1), \; t\in {\mathbb {R}}\\
v(x,0)=0,&x\in (0,1)\\
v(0,t)=v(1,t)=0,\quad &t\in{\mathbb {R}} \\
v_{xx}(0,t)=   v_{xx}(1,t)=0, &t\in{\mathbb {R}}
\end{array}\right.
$$
with
\begin{equation}\label{nonlinearterm}
f(x,t)=-i\Big ((1-x)(h_1'-\frac{1}{6} h_5')(t)+\frac{1}{6}(1-x)^3 h_5'(t)\Big ).
\end{equation}
By odd extension,  $f(x,t)$ can be expressed as
$$
f(x,t)=\sum_{k=1}^\infty \beta_k(t) \sin(k\pi x), \qquad\hbox{with}\qquad \beta_k(t)=2\int_0^1 f(x,t)\sin(k\pi x)dx.
$$
By \eqref{nonlinearterm},
\begin{eqnarray*}
	\beta_k(t) &=& -2i\int_0^1 (1-x)^3\sin(k\pi x) dx\frac{1}{6} h_5'(t)\\
	&=&-2i \int_0^1 (1-x)\sin(k\pi x) dx (h_1'-\frac{1}{6} h_5')(t) \\
	&=& -2i \Big (\frac{1}{k\pi}h_1'(t)- \frac{1}{(k\pi)^3}  h_5'(t)\Big ).
\end{eqnarray*}
Write $v(x,t)$ as
\begin{equation}\label{fomular for v}
v(x,t)=\sum_{k=1}^\infty\alpha_k(t)\sin k\pi x.
\end{equation}
Then, for $k=1,2,\cdots,$
\begin{equation}\label{fourier series}
i\frac{d}{dt}\alpha_k(t)+(k\pi)^4\alpha_k(t)=\beta_k(t), \qquad \alpha_k(0)=0.
\end{equation}
We have that
$$
\alpha_k(t)=2\int_0^te^{i(k\pi)^4(t-\tau)}\Big (\frac{1}{k\pi}h_1'(\tau)- \frac{1}{(k\pi)^3}  h_5'(\tau)\Big )d\tau.
$$
Taking into account that $h_1(0)=0$ and $ h_5(0)=0$, we have that
$$
\alpha_k(t)=2\Big (\frac{1}{k\pi}h_1(t)- \frac{1}{(k\pi)^3}  h_5(t)\Big )+\int_0^t e^{i(k\pi)^4(t-\tau)}2i(k\pi)^4\Big (\frac{1}{k\pi}h_1(\tau)- \frac{1}{(k\pi)^3}  h_5(\tau)\Big )d\tau.
$$
From the construction, we have
$$
\sum_{k=1}^\infty 2\Big (\frac{1}{k\pi}h_1(t)- \frac{1}{(k\pi)^3}  h_5(t)\Big )\sin(k\pi x) =-\Big ((1-x)(h_1-\frac{1}{6} h_5)(t)+\frac{1}{6}(1-x)^3 h_5(t)\Big ).
$$
Substituting $\alpha_k(t)$ into the original Fourier series representation  and taking above equation into account, it yields
$$	\begin{array}{lll}
v(x,t) &=& 
\displaystyle
-(1-x)(h_1-\frac{1}{6} h_5)(t)-\frac{1}{6}(1-x)^3 h_5(t)\\
&&
\displaystyle
+ \sum_{k=1}^\infty \int_0^t e^{i(k\pi)^4(t-\tau)}2i(k\pi)^4\Big (\frac{1}{k\pi}h_1(\tau)- \frac{1}{(k\pi)^3}  h_5(\tau)\Big )d\tau\sin(k\pi x)  \\
\end{array}
$$
which in turn implies that
$$
u(x,t)=\sum_{k=1}^\infty \int_0^t e^{i(k\pi)^4(t-\tau)}2i(k\pi)^4\Big (\frac{1}{k\pi}h_1(\tau)- \frac{1}{(k\pi)^3}  h_5(\tau)\Big )d\tau\sin(k\pi x).
$$
Next, consider the case of $h_1\equiv  h_5\equiv 0$ and $h_2(0)=h_6(0)=0$. Let $\tilde x=1-x$, we have the same situation as we just studied. Thus, if $h_1\equiv  h_5\equiv 0$ and $h_2(0)=h_6(0)=0$,
$$
u(x,t)=\sum_{k=1}^\infty (-1)^{k+1}\int_0^t e^{i(k\pi)^4(t-\tau)}2i(k\pi)^4\Big (\frac{1}{k\pi}h_2(\tau)- \frac{1}{(k\pi)^3} h_6(\tau)\Big )d\tau\sin(k\pi x).
$$
Combining the above two cases, due to the linearity of the system, we finish the proof of  Proposition \ref{expruxt-02}.

We now consider the boundary integral
$$	\begin{array}{lll}
u_{0,h}=W_{0,N} h &=& 
\displaystyle
\sum_{k=1}^\infty 2i(k\pi)^3\int_0^t e^{i(k\pi)^4(t-\tau)}h(\tau)d\tau  \sin(k\pi x)\\
&=& 
\displaystyle
\sum_{k=-\infty}^\infty (k\pi)^3e^{i(k\pi)^4t}\int_0^t e^{-i(k\pi)^4\tau}h(\tau)d\tau  e^{i k\pi x}.
\end{array}
$$
\begin{proposition}\label{p4.6-1}
	For any $h\in H^{\frac{3}{4}}(0,T)$, let $u=W_{0,N} h$ and $\Omega_T=(0,1)\times (0,T)$. Then $u$ belongs to $L^4(\Omega_T)\cap C([0,T]; L^2(0,1))$ and satisfies
	$$
	\norm{ u_{0,h}}_{ L^4(\Omega_T)}\leq C\norm{h}_{H^{\frac{3}{4}}(0,T)}
	$$
	and
	$$
	\sup_{0\leq t\leq T}\norm{ u_{0,h}(\cdot, t)}_{L^2(0,1)}\leq C\norm{h}_{H^{\frac{3}{4}}(0,T)}.
	$$
\end{proposition}

{\bf Proof:}  Without loss of generality, we assume that  $h(t)=0$ for $t\notin(0,T)$.   Let $h(\tau)=\int_{-\infty}^{\infty}e^{i\lambda\pi^4 \tau}\hat h(\lambda)d\lambda.$ Set $\alpha_k=i(\lambda-k^4)\pi^4,\beta_k=(k\pi)^3$. Then $u_{0,h}$ has the form:
\begin{eqnarray}
\nonumber    u_{0,h} & = & \displaystyle\sum_{k=-\infty}^\infty \beta_k e^{i(k\pi)^4t} e^{i k\pi x}\int_{ -\infty}^\infty\hat h(\lambda)\int_0^t e^{i\lambda\pi^4 \tau-i(k\pi)^4\tau}d\tau  \\
\nonumber     & = &\displaystyle \sum_{k=-\infty}^\infty \beta_k e^{i(k\pi)^4t} e^{i k\pi x}\int_{ -\infty}^\infty\hat h(\lambda) \frac{e^{\alpha_k t}-1}{\alpha_k}d\lambda \\
\nonumber    & = & \displaystyle\sum_{k=-\infty}^\infty \beta_k e^{i(k\pi)^4t} e^{i k\pi x}\Big (\int_{ -\infty}^0+\int_{ 0}^\infty\Big )\hat h(\lambda) \frac{e^{\alpha_k t}-1}{\alpha_k}d\lambda \\
\nonumber     & = & \displaystyle I^-(x,t)+I^+(x,t).
\end{eqnarray}
For $I^+(x,t)$, we have
\begin{eqnarray*}
	\nonumber I^+(x,t) &=&  \sum_{k=-\infty}^\infty \beta_k e^{i(k\pi)^4t} e^{i k\pi x} \int_{ 0}^\infty \hat h(\lambda)\psi(\lambda-k^4)\sum_{n=1}^\infty \frac{(\alpha_k t)^n}{ n! \alpha_k}d\lambda  \\
	\nonumber   & & +\sum_{k=-\infty}^\infty \beta_k  e^{i k\pi x} \int_{ 0}^\infty \hat h(\lambda) (1-\psi(\lambda-k^4))\frac{e^{i\lambda\pi^4 t}}{\alpha_k}d\lambda  \\
	\nonumber   & & -\sum_{k=-\infty}^\infty \beta_k  e^{i(k\pi)^4t} e^{i k\pi x} \int_{ 0}^\infty \hat h(\lambda)(1-\psi(\lambda-k^4)) \frac{1}{\alpha_k}d\lambda  \\
	\nonumber  &=& I_1^++ I_2^++I_3^+,
\end{eqnarray*}
where $\psi$ is a bump function associated to interval $[-4,4]$.  $I_1^+$ can be expressed as the form
$$
I_1^+=\sum_{n=1}^\infty\frac{(-i\pi^4)^{n}}{n!} I_{1,n}^+t^{n}
$$
with
\begin{equation}\label{i1+}
I_{1,n}^+=\sum_{k=-\infty}^\infty \beta_k e^{i(k\pi)^4t} e^{i k\pi x} \int_{ 0}^\infty \hat h(\lambda)\psi(\lambda-k^4)(k^4-\lambda)^{n-1} d\lambda.
\end{equation}
Using Proposition \ref{L4estimates}, we get
\begin{eqnarray*}
	\norm{I_{1,n}^+}^2_{L^4(\Omega_T)\cap L^\infty(0,T; L^2(0,1))}&\leq&C\Big (\sum_{k=-\infty}^\infty|\beta_k|^2\|\int_0^\infty \hat h(\lambda)\psi(k^4-\lambda)(k^4-\lambda)^{n-1}d\lambda\|^2\Big )\\
	&\leq& C 4^n\Big (\sum_{k=-\infty}^\infty k^{6}\|\int_{|\lambda-k^4|\leq 4}\hat h(\lambda)d\lambda\|^2\Big )\\
	&\leq& C 4^n\Big (\sum_{k=-\infty}^\infty k^{6}\int_{|\lambda-k^4|\leq 4}|\hat h(\lambda)|^2d\lambda\Big )\\
	&\leq& C4^n\Big (\sum_{k=-\infty}^\infty \int_{|\lambda-k^4|\leq 4}||
	\lambda|^{\frac{3}{4}}\hat h(\lambda)|^2d\lambda\Big )\\
	&\leq& C 4^n\norm{h}^2_{H^{\frac{3}{4}}({\mathbb {R}})}.
\end{eqnarray*}
Taking above inequality into \eqref{i1+} with $0\leq t\leq T$, we have the estimate of $I_1^+$.

Now we estimate  $I_2^+$.  Rewrite $I_2^+$ as
\begin{eqnarray*}
	&&I_2^+=\sum_{k=1}^\infty 2\beta_k\sin(k\pi x)\int_{ 0}^\infty \hat h(\lambda) (1-\psi(\lambda-k^4))\frac{e^{i\lambda\pi^4 t}}{(\lambda-k^4)\pi^4}d\lambda \\
	&=& \sum_{k=1}^\infty \frac{\beta_k}{k^2\pi^4}\sin(k\pi x)\int_{ 0}^\infty \hat h(\lambda) (1-\psi(\lambda-k^4)) e^{i\lambda\pi^4 t}\Big (\frac{1}{\sqrt\lambda-k^2}-\frac{1}{\sqrt\lambda+k^2}\Big )d\lambda  \\
	&=& \sum_{k=1}^\infty \frac{\beta_k}{2k^3\pi^4}\sin(k\pi x)\int_{ 0}^\infty \hat h(\lambda) (1-\psi(\lambda-k^4)) e^{i\lambda\pi^4 t}\Big (\frac{1}{\sqrt[4]\lambda-k}-\frac{1}{\sqrt[4]\lambda+k}-\frac{2k}{\sqrt\lambda+k^2}\Big )d\lambda  \\
	&=& \sum_{k=1}^\infty \frac{2\beta_k}{k^3\pi^4}\sin(k\pi x)\int_{ 0}^\infty \mu^3\hat h(\mu^4) (1-\psi(\mu^4-k^4)) e^{i\mu^4\pi^4 t}\Big (\frac{1}{\mu-k}-\frac{1}{\mu+k}-\frac{2k}{\mu^2+k^2}\Big )d\mu.  \\
\end{eqnarray*}
It is easy to see the main term is the $\frac{1}{\mu-k}$ term. Applying Lemma \ref{A-1-1} with $f(\mu)=\mu^3\hat h(\mu^4)$ in Appendix, it holds
\begin{eqnarray*}
	&&\sup_{0\leq t\leq T}\norm{I_2^+(\cdot, t)}_{L^2(0,1)}^2\\
	&\leq & C\sum_{k=1}^\infty \|\int_{ 0}^\infty \mu^3\hat h(\mu^4) (1-\psi(\mu^4-k^4)) \Big (\frac{1}{\mu-k}-\frac{1}{\mu+k}-\frac{2k}{\mu^2+k^2}\Big )d\mu\|^2\\
	&\leq&C\norm{(|\mu|+1)^{\frac{3}{4}}\hat h(\mu) }^2_{({\mathbb {R}}^+)}\leq \norm{h}^2_{H^{\frac{3}{4}}({\mathbb {R}}^+)}.
\end{eqnarray*}

To estimate the $L^4(\Omega_T)-$norm, we write $I_2^+(x,t)$ as
\begin{eqnarray}\label{i212}
I_2^+&=&  \sum_{k=-\infty}^\infty \beta_k  e^{i k\pi x} \Big (\int_{ 0}^{\frac{k^4}{2}} +\int_{\frac{k^4}{2}}^\infty \Big )\hat h(\lambda) (1-\psi(\lambda-k^4))\frac{e^{i\lambda\pi^4 t}}{\alpha_k}d\lambda  \\
\nonumber  &=& I_{2,1}^++ I_{2,2}^+.
\end{eqnarray}
Applying Lemma \ref{6.5} in Appendix, by taking $\varepsilon>0$ sufficiently small, the second term in \eqref{i212} satisfies
\begin{eqnarray*}
	&&\nonumber \norm{I_{2,2}^+}_{L^4(\Omega_T)}^2\\
	\nonumber
	&\leq& C\Big (  \sum_{k=-\infty}^\infty\int_{{\mathbb {R}}}|\beta_k|^2  \frac{|\hat h(\lambda)|^2}{(|\lambda-k^4|+1)^2} (|\lambda-k^4|+1)^{\frac{1}{2}+\varepsilon}\chi_{[\frac{k^4}{2},\infty)}(\lambda)(1-\psi(\lambda-k^4))^2 d\lambda\Big )  \\
	\nonumber  &\leq&C  \int_0^\infty |\lambda|^{\frac{3}{2}}  |\hat h(\lambda)|^2 \sum_{k=-\infty}^\infty  \frac{1}{(|\lambda-k^4|+1)^{\frac{3}{2}-\varepsilon} } d\lambda \leq C  \int_0^\infty |\lambda|^{\frac{3}{2}}  |\hat h(\lambda)|^2 d\lambda\leq C\norm{h}_{H^{\frac{3}{4}}({\mathbb {R}})}^2.
\end{eqnarray*}
Moreover, 
\begin{eqnarray*}
	\nonumber
	&&\|I_{2,1}^+\|\\
	\nonumber&=& \|2\int_0^\infty\Big (\sum_{k=1}^\infty \chi_{[0, \frac{k^4}{2}]}(\lambda)(1-\psi(k^4-\lambda))\frac{|\beta_k|\sin(k\pi x)}{|\alpha_k|}\Big ) e^{i\lambda\pi^4 t} \hat h(\lambda)d\lambda\|  \\
	\nonumber
	&\leq& \|\frac{1}{\pi}\int_0^\infty e^{i\lambda\pi^4 t} \hat h(\lambda)\Big (\sum_{k=1}^\infty \sin(k\pi x)\chi_{[0, \frac{k^4}{2}]}(\lambda)(1-\psi(k^4-\lambda))\Big (\frac{1}{\sqrt[4]{\lambda}-k}-\frac{1}{\sqrt[4]{\lambda}+k}-\frac{2k}{\sqrt{\lambda}+k^2}\Big )d\lambda\| \\
	\nonumber
	&=& \|\frac{1}{\pi}\int_0^\infty e^{i\lambda\pi^4 t} \hat h(\lambda)\Big (\sum_{k=\lfloor\sqrt[4]{2\lambda}\rfloor}^\infty \sin(k\pi x)\chi_{[0, \frac{k^4}{2}]}(\lambda)\Big (\frac{1}{\sqrt[4]{\lambda}-k}-\frac{1}{\sqrt[4]{\lambda}+k}-\frac{2k}{\sqrt{\lambda}+k^2}\Big )d\lambda\| \\
	\nonumber
	&\leq&\frac{1}{\pi}\int_0^\infty |\hat h(\lambda)|\|\sum_{k=\lfloor\sqrt[4]{2\lambda}\rfloor}^\infty \sin(k\pi x)\Big (\frac{1}{\sqrt[4]{\lambda}-k}-\frac{1}{\sqrt[4]{\lambda}+k}-\frac{2k}{\sqrt{\lambda}+k^2}\Big )\| d\lambda.
\end{eqnarray*}
Use Lemma \ref{p37-middle} in Appendix for the first term  and similar arguments for the others, for any $\bar\alpha\in (\frac{1+\alpha}{4}, \frac{3}{4}]$
\footnote%
{
	It is sufficient that $\frac{2-2\alpha}{4}+2\bar \alpha >1\Leftrightarrow \bar \alpha>\frac{1+\alpha}{4}$.
}
we have
$$
\begin{array}{lll}
\displaystyle\|I_{2,1}^+(x,t)\|
&\leq &
\displaystyle C|x|^{\alpha-1}\int_0^\infty \frac{|\hat h(\lambda)|}{(1+\sqrt[4]{\lambda})^{1-\alpha}}d\lambda\\
&\leq &
\displaystyle C|x|^{\alpha-1}\int_0^\infty (1+|\lambda|)^{\bar\alpha}\frac{|\hat h(\lambda)|}{(1+\sqrt[4]{\lambda})^{1-\alpha}(1+|\lambda|)^{\bar\alpha}}d\lambda \\
&\leq &
\displaystyle C|x|^{\alpha-1}\Big (\int_0^\infty (1+|\lambda|)^{2\bar\alpha} |\hat h(\lambda)|^2 d\lambda\Big )^{\frac{1}{2}}
\Big (\int_0^\infty  \frac{d\lambda }{(1+\sqrt[4]{\lambda})^{2-2\alpha}(1+|\lambda|)^{2\bar\alpha}}\Big )^{\frac{1}{2}}\\
&\leq &
\displaystyle C|x|^{\alpha-1}\norm{h}_{H^{\bar\alpha}(\mathbb{R})}.
\end{array}
$$
The above estimates leads to
$$
\norm{I_2^+}^2_{L^4(\Omega_T)}\leq C\norm{h}^2_{H^{\frac{3}{4}}(\mathbb{R})}.
$$

Now we study $I_3^+(x,t)$. We have
$$
\begin{array}{lll}
&&\displaystyle\norm{I_{3}^+}^2_{L^4(\Omega_T)\cap L^\infty(0,T; L^2(0,1))}\\
&\leq &
\displaystyle C \Big (\sum_{k=1}^\infty k^6\|\int_0^\infty \hat h(\lambda)\frac{1-\psi(k^4-\lambda)}{\lambda-k^4}d\lambda\|^2\Big )\\
&\leq &
\displaystyle C \sum_{k=1}^\infty \|\int_0^\infty \hat h(\lambda)\Big (\frac{1}{\sqrt[4]{\lambda}-k}-\frac{1}{\sqrt[4]{\lambda}+k}-\frac{2k}{\sqrt{\lambda}+k^2}\Big )(1-\psi(k^4-\lambda))d\lambda\|^2\\
&\leq &
\displaystyle C \Big (\sum_{k=1}^\infty \|\int_0^\infty \mu^3\hat h(\mu^4)\frac{1}{\mu-k}(1-\psi(k^4-\mu^4))d\lambda\|^2\\
&&
\displaystyle\quad+\sum_{k=1}^\infty \|\int_0^\infty \mu^3\hat h(\mu^4)\frac{1}{\mu+k}(1-\psi(k^4-\mu^4))d\lambda\|^2\\
&&
\displaystyle\quad+ \sum_{k=1}^\infty \|\int_0^\infty \mu^3\hat h(\mu^4)\frac{2k}{\mu^2+k^2}(1-\psi(k^4-\mu^4))d\lambda\|^2\Big ).
\end{array}
$$
Use Lemma \ref{A-1-1} with $f(\mu)=\mu^3 \hat h(\mu^4)$ and Lemma \ref{p2.1},  it holds
$$
\begin{array}{lll}
& & \displaystyle
\hspace{5mm}\norm{I_{3}^+}^2_{L^4(\Omega_T)\cap L^\infty(0,T; L^2(0,1))}\\
&&\leq 
\displaystyle C \left\{\int_0^\infty (\mu+1)^3\mu^6|\hat h(\mu^4)|^2d\mu
\right.\\
& &
\left.\displaystyle
\hspace{5mm}
+\int_0^\infty \|\Big (\int_0^\infty \frac{\mu^3\hat h(\mu^4)}{\mu+y}d\mu\Big )^2+\Big (\int_0^\infty \frac{\mu^3\hat h(\mu^4)}{\mu^2+y^2}d\mu\Big )^2\|dy\right\}\\
&& \leq C\norm{h}^2_{H^{\frac{3}{4}}(\mathbb{R})}.
\end{array}
$$
In sum, it appears that
$$
\norm{I^+}^2_{L^4}\leq C\norm{h}^2_{H^{\frac{3}{4}}(\mathbb{R})}.
$$

The estimates for $I^-(x,t)$ follows in the same way and notice that in this case there is no case $|\lambda-k^4|<4$ for $k\neq0,1,-1$ since $\lambda<0$.
We finish the proof.
\vskip 5mm

Actually, we can get some smoother estimates.

\begin{proposition}\label{p4.6-2}
	Let $0\leq s\leq4$ be given. For any $h\in H_0^{\frac{3+s}{4}}(0,T)$, let $u=W_{0,N} h$. Then $\partial_x^s u$ belongs to $L^4(\Omega_T)\cap C([0,T]; L^2(0,1))$ and satisfies
	$$
	\norm{\partial_x^s u}_{ L^4(\Omega_T)}\leq C\norm{h}_{H^{\frac{3+s}{4}}(0,T)}
	$$
	and
	$$
	\sup_{0\leq t\leq T}\norm{\partial_x^s u(\cdot, t)}_{L^2(0,1)}\leq C\norm{h}_{H^{\frac{3+s}{4}}(0,T)}.
	$$
\end{proposition}
\begin{remark} Note that the $H^{\frac{3+s}{4}}_0(0,T)$ denotes the interpolation space between $H^{\frac{3}{4}}$ and $H^\frac{7}{4}$. The condition $s\leq 4$ is just a technical requirement. We do not pursue for more smooth condition. It is enough for our argument.
\end{remark}

{\bf Proof:} We only give the proof for $s=4$. $0\leq s\leq 4$ follows from the interpolation argument. We here reform our goal to set up
$$
\norm{\partial_x^4 u}_{ L^4(\Omega_T)}\leq C\norm{h}_{H^{\frac{7}{4}}(0,T)}
$$
and
$$
\sup_{0\leq t\leq T}\norm{\partial_x^4 u(\cdot, t)}_{L^2(0,1)}\leq C\norm{h}_{H^{\frac{7}{4}}(0,T)}.
$$

From \eqref{u-vrelation}, and the definition of $f$ with $ h_5\equiv0$, we always have
$$\partial_x^4u=\partial_x^4v.$$
Noticing \eqref{fomular for v} and \eqref{fourier series}, we have
$$
\norm{\partial_x^4u}_{ L^4(\Omega_T)\bigcap L^\infty(0,T; L^2(0,1))}\leq \norm{\partial_t v}_{ L^4(\Omega_T)\bigcap L^\infty(0,T; L^2(0,1))}+\norm{f}_{L^4(\Omega_T)\bigcap L^\infty(0,T; L^2(0,1))}.
$$
By the definition of $f$, and Sobolev embedding
$$\norm{f}_{ L^4(\Omega_T)}\leq C\norm{h}_{L^4(0,T)}\leq C\norm{h}_{H^\frac{7}{4}(0,T)}$$
and
$$\norm{f}_{L^\infty((0,T)L^2(0,1))}\leq C\norm{h}_{L^\infty(0,T)}\leq C\norm{h}_{H^\frac{7}{4}(0,T)}.$$
On the other hand, going through the argument of Proposition \ref{p4.6-1} again, we can obtain
$$\norm{\partial_t v}_{L^4(\Omega_T)\bigcap L^\infty(0,T; L^2(0,1))}\leq C\norm{h'}_{H^{\frac{3}{4}}(0,T)}=C\norm{h}_{H^{\frac{7}{4}}(0,T)}.$$
We finish our proof.
\vskip 5mm

Now we consider the second order boundary integral
$$
u_{2,h}=W_{2,N} h= \sum_{k=1}^\infty (-2ik\pi)\int_0^t e^{i(k\pi)^4(t-\tau)}h(\tau)d\tau \sin(k\pi x).
$$ 
Again, by the same argument, we can obtain the following estimates for $u_{2,h}.$

\begin{proposition}\label{p4.6-3}
	For any given $T>0$. Let $\Omega_T=(0,1)\times(0,T)$. If $h\in H^{\frac{1}{4}}(0,T)$, then we have
	$$
	u_{2,h}=W_{2,N}(\cdot)h\in L^4(\Omega_T)\cap C([0,T]; L^2(0,1))
	$$
	and there is a positive constant $C_T$ depending only on $T$ such that
	$$
	\norm{u_{2,h}}_{ L^4(\Omega_T)}\leq C_T\norm{h}_{H^{\frac{1}{4}}(0,T)}
	$$
	and
	$$
	\sup_{0\leq t\leq T}\norm{u_{2,h}(\cdot, t)}_{L^2(0,1)}\leq C_T\norm{h}_{H^{\frac{1}{4}}(0,T)}.
	$$
\end{proposition}
{\bf Proof:} Set $\alpha_k=i(\lambda-k^4)\pi^4, \beta_k=-2k\pi$. Similar to Proposition \ref{p4.6-1},  $u_{2,h}$ has the form
$$
u_{2,h}=I^{-}(x,t)+I_1^{+}(x,t)+I_2^{+}(x,t)+I_3^{+}(x,t)
$$
with
\begin{eqnarray}
\nonumber
I^{-}(x,t) &=&\sum_{k=-\infty}^\infty \beta_k e^{i(k\pi)^4t} e^{i k\pi x}\int_{ -\infty}^0\hat h(\lambda) \frac{e^{\alpha_k t}-1}{\alpha_k}d\lambda,\\
\nonumber
I_1^{+}(x,t)&=&\sum_{k=-\infty}^\infty \beta_k e^{i(k\pi)^4t} e^{i k\pi x} \int_{ 0}^\infty \hat h(\lambda)\psi(\lambda-k^4)\sum_{n=1}^\infty \frac{(\alpha_k t)^n}{ n! \alpha_k}d\lambda,  \\
\nonumber
I_2^{+}(x,t)&=& \sum_{k=-\infty}^\infty \beta_k  e^{i k\pi x} \int_{ 0}^\infty \hat h(\lambda) (1-\psi(\lambda-k^4))\frac{e^{i\lambda\pi^4 t}}{\alpha_k}d\lambda,  \\
\nonumber
I_3^{+}(x,t)&=& -\sum_{k=-\infty}^\infty \beta_k  e^{i(k\pi)^4t} e^{i k\pi x} \int_{ 0}^\infty \hat h(\lambda)(1-\psi(\lambda-k^4)) \frac{1}{\alpha_k}d\lambda.
\end{eqnarray}

Recalling \eqref{i1+},  we have the bound of $I_1^{+}(x,t)$ due to the fact that
\begin{eqnarray*}
	\norm{I_{1,n}^+}^2_{L^4(\Omega_T)\cap L^\infty(0,T; L^2(0,1))}&\leq&C\Big (\sum_{k=-\infty}^\infty|\beta_k|^2\|\int_0^\infty \hat h(\lambda)\psi(k^4-\lambda)(k^4-\lambda)^{n-1}d\lambda\|^2\Big )\\
	&\leq& C  4^n \Big (\sum_{k=-\infty}^\infty k^{2}\|\int_{|\lambda-k^4|\leq 4}\hat h(\lambda)d\lambda\|^2\Big )\\
	&\leq& C 4^n \Big (\sum_{k=-\infty}^\infty k^{2}\int_{|\lambda-k^4|\leq 4}|\hat h(\lambda)|^2d\lambda\Big )\\
	&\leq& C 4^n \Big (\sum_{k=-\infty}^\infty \int_{|\lambda-k^4|\leq 4}||
	\lambda|^{\frac{1}{4}}\hat h(\lambda)|^2d\lambda\Big )\\
	&\leq& C 4^n  \norm{h}^2_{H^{\frac{1}{4}}({\mathbb {R}})}.
\end{eqnarray*}
For $I_2^{+}(x,t)$, it holds
$$
I_2^+=\sum_{k=1}^\infty \frac{2\beta_k}{k\pi^4}\sin(k\pi x)\int_{ 0}^\infty \frac{\mu^3}{k^2}\hat h(\mu^4) (1-\psi(\mu^4-k^4)) e^{i\mu^4\pi^4 t}\Big (\frac{1}{\mu-k}-\frac{1}{\mu+k}-\frac{2k}{\mu^2+k^2}\Big )d\mu.
$$
Similar computation shows that
$$
\sup_{0\leq t\leq T}\norm{I_2^+(\cdot, t)}_{L^2(0,1)}^2
\leq C\norm{(|\mu|+1)^{\frac{1}{4}}\hat h(\mu) }^2_{L^2(0,1)}\leq \norm{h}^2_{H^{\frac{1}{4}}({\mathbb {R}}^+)}.
$$
The estimate of its $L^4$-norm is similar with the notations $I_{2,2}^+$ and $I_{2,2}^+$ (see \eqref{i212}).  In fact,  applying  Lemma \ref{6.5}, it holds for sufficient small $\varepsilon>0$,
\begin{eqnarray*}
	\nonumber \norm{I_{2,2}^+}_{L^4(\Omega_T)}^2&\leq& C\Big (  \sum_{k=-\infty}^\infty\int_{{\mathbb {R}}}|\beta_k|^2  \frac{|\hat h(\lambda)|^2}{(|\lambda-k^4|+1)^2} (|\lambda-k^4|+1)^{\frac{1}{2}+\varepsilon}\chi_{[\frac{k^4}{2},\infty)}(\lambda)(1-\psi(\lambda-k^4))^2 d\lambda\Big )  \\
	\nonumber  &\leq&C  \int_0^\infty |\lambda|^{\frac{1}{2}}  |\hat h(\lambda)|^2 \sum_{k=-\infty}^\infty  \frac{1}{(|\lambda-k^4|+1)^{\frac{1}{2}-\varepsilon} } d\lambda \leq C  \int_0^\infty |\lambda|^{\frac{1}{2}}  |\hat h(\lambda)|^2 d\lambda\leq C\norm{h}_{H^{\frac{1}{4}}({\mathbb {R}})}^2.
\end{eqnarray*}

For $|I_{2,1}^+|$,
\begin{eqnarray*}
	\nonumber
	\|I_{2,1}^+\|&=& \|2\int_0^\infty\Big (\sum_{k=1}^\infty \chi_{[0, \frac{k^4}{2}]}(\lambda)(1-\psi(k^4-\lambda))\frac{|\beta_k|\sin(k\pi x)}{|\alpha_k|}\Big ) e^{i\lambda\pi^4 t} \hat h(\lambda)d\lambda\|  \\
	\nonumber
	&\leq& \|\frac{1}{\pi}\int_0^\infty e^{i\lambda\pi^4 t} \hat h(\lambda)\Big (\sum_{k=1}^\infty \frac{\sin(k\pi x)}{k^2}\chi_{[0, \frac{k^4}{2}]}(\lambda)(1-\psi(k^4-\lambda)\Big ) \\
	\nonumber
	&&\hspace{38mm}
	\Big (\frac{1}{\sqrt[4]{\lambda}-k}-\frac{1}{\sqrt[4]{\lambda}+k}-\frac{2k}{\sqrt{\lambda}+k^2}\Big )d\lambda\| \\
	\nonumber
	&=& \|\frac{1}{\pi}\int_0^\infty e^{i\lambda\pi^4 t} \hat h(\lambda)\Big (\sum_{k=\lfloor\sqrt[4]{2\lambda}\rfloor}^\infty \frac{\sin(k\pi x)}{k^2}\chi_{[0, \frac{k^4}{2}]}(\lambda)\Big (\frac{1}{\sqrt[4]{\lambda}-k}-\frac{1}{\sqrt[4]{\lambda}+k}-\frac{2k}{\sqrt{\lambda}+k^2}\Big )d\lambda\| \\
	\nonumber
	&\leq&\frac{1}{\pi}\int_0^\infty |\hat h(\lambda)|\|\sum_{k=\lfloor\sqrt[4]{2\lambda}\rfloor}^\infty \dfrac{\sin(k\pi x)}{k^2}\Big (\frac{1}{\sqrt[4]{\lambda}-k}-\frac{1}{\sqrt[4]{\lambda}+k}-\frac{2k}{\sqrt{\lambda}+k^2}\Big )\| d\lambda.
\end{eqnarray*}
Use Lemma \ref{p37-middle} in Appendix for the first term  and similar arguments for the others, for any $\bar\alpha\in (\frac{1+\alpha}{4}, \frac{1}{4}]$  
we have
$$
\|I_{2,1}^+(x,t)\|
\leq  C|x|^{\alpha-1}\norm{h}_{H^{\bar\alpha}(\mathbb{R}^+)}.
$$
and leads to
$$
\norm{I_2^+}^2_{L^4(\Omega_T)}\leq C\norm{h}^2_{H^{\frac{1}{4}}(\mathbb{R}^+)}.
$$

For $I_3^+(x,t)$, we have
$$
\begin{array}{lll}
\displaystyle&&\norm{I_{3}^+}^2_{L^4(\Omega_T)\cap L^\infty(0,T; L^2(0,1))}
\leq
\displaystyle C \Big (\sum_{k=1}^\infty k^2\|\int_0^\infty \hat h(\lambda)\frac{1-\psi(k^4-\lambda)}{\lambda-k^4}d\lambda\|^2\Big )\\
&\leq &
\displaystyle C \sum_{k=1}^\infty \|\int_0^\infty \hat h(\lambda)\Big (\frac{1}{\sqrt{\lambda}-k^2}-\frac{1}{\sqrt{\lambda}+k^2}\Big )(1-\psi(k^4-\lambda))d\lambda\|^2\\
&\leq &
\displaystyle C \Big (\sum_{k=1}^\infty \|\int_0^\infty \mu^3\hat h(\mu^4)\frac{1}{\mu^2-k^2}(1-\psi(k^4-\mu^4))d\lambda\|^2\\
&&
\displaystyle\quad+ \sum_{k=1}^\infty \|\int_0^\infty \mu^3\hat h(\mu^4)\frac{1}{\mu^2+k^2}(1-\psi(k^4-\mu^4))d\lambda\|^2\Big ).
\end{array}
$$
Use Lemma \ref{A-1-1} with $f(\mu)=\mu^3 \hat h(\mu^4)$ and Lemma \ref{p2.1},  it holds
$$
\begin{array}{lll}
\displaystyle\norm{I_{3}^+}^2_{L^4(\Omega_T)\cap L^\infty(0,T; L^2(0,1))}
&\leq &
\displaystyle C \left\{\int_0^\infty (\mu+1)^3\mu^6|\hat h(\mu^4)|^2d\mu\right.\\
& &
\left.\displaystyle+\int_0^\infty \|\Big (\int_0^\infty \frac{\mu^3\hat h(\mu^4)}{\mu+y}d\mu\Big )^2+\Big (\int_0^\infty \frac{\mu^3\hat h(\mu^4)}{\mu^2+y^2}d\mu\Big )^2\|dy\right\}\\
&\leq& C\norm{h}^2_{H^{\frac{3}{4}}(\mathbb{R}^+)}.
\end{array}
$$
In sum, it appears that
$$
\norm{I^+}^2_{L^4}\leq C\norm{h}^2_{H^{\frac{1}{4}}(\mathbb{R}^+)}.
$$

The estimates for $I^-(x,t)$ follows in the same way and notice that in this case there is no case $|\lambda-k^4|<4$ for $k\neq0,1,-1$ since $\lambda<0$.
We finish the proof.
\vskip 5mm

Similar to Proposition \ref{p4.6-2}, it holds
\begin{proposition}\label{p4.6-4}
	Let $0\leq s\leq 4$ be given. For any $h\in H_0^{\frac{1+s}{4}}(0,T)$, let $u=W_{2,N} h$. Then $\partial_x^s u$ belongs to $L^4(\Omega_T)\cap C([0,T]; L^2(0,1))$ and satisfies
	$$
	\norm{\partial_x^s u}_{ L^4(\Omega_T)}\leq C\norm{h}_{H^{\frac{1+s}{4}}(0,T)}
	$$
	and
	$$
	\sup_{0\leq t\leq T}\norm{\partial_x^s u(\cdot, t)}_{L^2(0,1)}\leq C\norm{h}_{H^{\frac{1+s}{4}}(0,T)}.
	$$
\end{proposition}

\subsection{Case 2: With Dirichlet boundary conditions}\label{3.12}
We now consider the IBVP
\begin{equation}\label{LinearIBVPD}
\left\{\begin{array}{ll}
i\partial_t u+  \partial_x^4u=0, & x\in (0,1), \; t\in {\mathbb {R}}\\
u(x,0)=\varphi(x),&x\in (0,1)\\
u(0,t)=r_1(t), u(1,t)=r_2(t),\quad &t\in{\mathbb {R}} \\
\partial_x u(0,t)=r_3(t),  \partial_x u(1,t)=r_4(t), &t\in{\mathbb {R}}
\end{array}\right.
\end{equation}
for the linear Schr\"odinger equation . Here $r_i(t), i=1,2,3,4$ are given in Lemma \ref{lema1} and  \ref{lema2}, accordingly. In fact, by the well-considered design of the boundary terms, system \eqref{LinearIBVPD} is divided into two new systems
\begin{equation}\label{LinearIBVPDo}
\left\{\begin{array}{ll}
i\partial_tu_o(x,t)+\partial_x^4 u_o(x,t)=0, & x\in (0,1), \; t\in {\mathbb {R}}\\
u_o(x,0)=\varphi_o(x),&x\in (0,1)\\
u_o(0,t)=u_o(1,t)=0,\quad &t\in{\mathbb {R}} \\
\partial_x u_o(0,t)=r_3(t),  \partial_x u_o(1,t)=r_4(t), &t\in{\mathbb {R}}
\end{array}\right.
\end{equation}
and
\begin{equation}\label{LinearIBVPDe}
\left\{\begin{array}{ll}
i\partial_tu_e(x,t)+\partial_x^4 u_e(x,t)=0, & x\in (0,1), \; t\in {\mathbb {R}}\\
u_e(x,0)=\varphi_e(x),&x\in (0,1)\\
u_e(0,t)=r_1(t), u_e(1,t)=r_2(t),\quad &t\in{\mathbb {R}} \\
\partial_x u_{e}(0,t)=0,  \partial u_{e}(1,t)=0, &t\in{\mathbb {R}}
\end{array}\right.
\end{equation}
where $\varphi_e(x)$ and $\varphi_o(x)$ are given by \eqref{phioe}.  $u_0$ and $u_e$ is defined in \eqref{uoext}. Let  $$W(t)\varphi:=u(x,t)=u_o(x,t)+u_e(x,t)\quad (x,t)\in(0,1)\times[0,1]$$ It solves \eqref{LinearIBVPD} for $t\in[0,1].$
\begin{proposition}\label{L4estimatesD}
	Let $s\geq 0$. For any $\varphi\in H^s(0,1)$, the solution of \eqref{LinearIBVPD}  for $t\in[0,1]$  satisfies
	$$
	\norm{u}_{ C([0,1]; H^s(0,1))}\leq C\norm{\varphi}_{H^s(0,1)}.
	$$
\end{proposition}
{\bf Proof}: It follows directly form the Plancherel theorem.

We now consider the  corresponding non-homogeneous problem with vanished Dirichlet boundary. 
$$
\left\{\begin{array}{ll}
i\partial_t u+  \partial_x^4u=f, & x\in (0,1), \; t\in {\mathbb {R}}\\
u(x,0)=0,&x\in (0,1)\\
u(0,t)=u(1,t)=0,\quad &t\in{\mathbb {R}} \\
\partial_x u(0,t)=  \partial_x u(1,t)=0, &t\in{\mathbb {R}}
\end{array}\right.
$$
By Duhamel's principle, the solution is
$$
u(t)=-i\int_0^tW^D(t-\tau)f(\cdot, \tau)d\tau.
$$
Here $W^D(t)$ is the $C_0$-group in $ L^2(0,1)$ generated by the operator $A^Dv=i v''''$ with domain $\mathcal{D}(A^D)=\{u\in H^4(0,1): u(x)=u_{x}(x)=0, \; x=0, 1\} $.

\begin{proposition}\label{p4.6} 
	Let $s\in[0,4]$ and $T>0$ be given. Let
	$$
	u(t)=\int_0^t W^D(t-\tau)f(\cdot, \tau)d\tau
	$$
	with $f\in L^1(0,T; H^s(0,1))$  satisfying
	$$
	f(0,t)=f(1, t)=f_{x}(0,t)=f_{x}(1,t)\equiv0
	$$
	when $s>\frac{5}{2}$ and
	$$
	f(0,t)=f(1, t)\equiv0
	$$
	with $s\in (\frac{1}{2},\frac{5}{2})$. Then $u\in C([0,T]; H^s(0,1))$ and
	\begin{equation}
	\label{TPex2D}\norm{v}_{C([0,T]; H^s(0,1))}\leq  C_{T,s}\norm{f}_{L^1(0,T; H^s(0,1))}
	\end{equation}
	where the constant $C_{T,s}$ depends only on $s$ and $T$.
\end{proposition}
{\bf Proof: }The existence of the solution $u$ follows from the semigroup theory. For any fixed $t\in[0,T]$, the inequality \eqref{TPex2D}  follows from the Plancherel theorem if $s = 0$ and $s = 4$. When $0 < s < 4$,  these inequality is deduced from interpolation theory using the results for $s = 0$ and $s = 4$.

\vskip 3mm

Now we consider the linear problem with non-homogeneous Dirichlet boundary conditions, i.e.
\begin{equation}\label{4thNBD}
\left\{\begin{array}{ll}
i\partial_t u+  \partial_x^4u=0, & (x, t)\in (0,1)\times \mathbb{R}\\
u(x,0)=0,&x\in (0,1)\\
u(0,t)=h_1(t),\quad u(1,t)=h_2(t),\quad &t\in\mathbb{R} \\
\partial_x u(0,t)= h_3(t), \quad  \partial_x u(1,t)=h_4(t), &t\in\mathbb{R}.
\end{array}\right.
\end{equation}
We put the compatibility conditions
\begin{equation}\label{compa}
h_1(0)=h_2(0)=0, \qquad  h_3(0)=h_4(0)=0
\end{equation}
if it is necessary.

\begin{proposition}\label{BoIn4}
	The solution of \eqref{4thNBD} can be expressed as 
	$$
	u(x,t)=W_{0,D} h_1+W_{1,D}  h_3+\Big(W_{0,D}h_2+W_{1,D} h_4\Big)\left| _{x\rightarrow1-x}\right. 
	$$
	in which the boundary integral operator $W_{0,D}$ and $W_{1,D}$ are given by
	\begin{eqnarray}
	W_{0,D} h_1&=&
	\displaystyle\label{beta01}
	\sum_{k=1}^{\infty}\int_0^te^{i(k\pi)^4(t-\tau)}
	\beta_{01,k} h_1(\tau)
	d\tau\sin k\pi x\\
	&&\displaystyle 
	+\sum_{k=1}^{\infty}\int_0^te^{i(k\pi)^4(t-\tau)}
	\beta_{02,k}h_1(\tau)
	d\tau\cos k\pi x\\
	W_{1,D}  h_3 &=&\displaystyle
	\sum_{k=1}^{\infty}\int_0^te^{i(k\pi)^4(t-\tau)}
	\beta_{11,k}h_3(\tau)
	d\tau\sin k\pi x\\
	&&\displaystyle 
	+\sum_{k=1}^{\infty}\int_0^te^{i(k\pi)^4(t-\tau)}
	\beta_{12,k} h_3(\tau)
	d\tau\cos k\pi x.
	\end{eqnarray}
	with
	\begin{equation}\label{betaij}
	\begin{array}{lll}
	&\hspace{-3mm} 
	\beta_{01,k}=-i(k\pi)^3-6ik\pi(\cos k\pi+1),
	&\beta_{02,k}=12i(k\pi-1),\\
	&\hspace{-3mm} 
	\beta_{11,k}=-2ik\pi(\cos k\pi+2),
	&\beta_{12,k}=i(k\pi)^2+6i(\cos k\pi-1).\\
	\end{array}
	\end{equation}
\end{proposition}
{\bf Proof :}  Similar to the case with Navier boundary conditions, we  need only to prove that  
\begin{equation}\label{uww}
u(x,t)=  W_{0,D} h_1+W_{1,D}  h_3
\end{equation}
under the assumption that $h_2(t)\equiv h_4(t)\equiv 0$ and $h_1(0)=h_3(0)=0$. The proof is divided into two steps.

{\bf Step 1:}   Define $v$ by
\begin{equation}\label{u-v}
v(x,t)=u(x,t)-f(x,t),
\end{equation}
with
\begin{equation}\label{fxt1}
f(x,t)= (1-x)^2(3h_1(t)+h_3(t))-(1-x)^3(2h_1(t)+h_3(t)).
\end{equation}
By direct computation, we need to solve
\begin{equation}\label{4th_Boundary_data_v-02}
\left\{\begin{array}{ll}
i\partial_t v+ \partial_x^4 v=-i\partial_t f
& x\in (0,1), \; t\in {\mathbb {R}}\\
v(x,0)=0,&x\in (0,1)\\
v(0,t)=v(1,t)=0,\quad &t\in{\mathbb {R}} \\
v_{x}(0,t)=   v_{x}(1,t)=0, &t\in{\mathbb {R}}.
\end{array}\right.
\end{equation}
For fixed $t$, we have $f\in C^\infty[0,1]$. By the same extension as in \eqref{phioe}, and modifying the value of $f(x,t)$ at $x=-1,0,1$ if necessary, we have	\begin{equation}\label{exfxt1}
f(x,t)=p_0(t)+\sum_{k=1}^{\infty}p_k(t)\cos k\pi x+\sum_{k=1}^{\infty}q_k(t)\sin k\pi x,
\end{equation}
with
$$
p_0(t)=\frac{1}{2}h_1-\frac{1}{12}h_3,
\quad
q_k(t)=b_{01,k}h_1+b_{11,k}h_3,
\quad
p_k(t)=b_{02,k}h_1+b_{12,k}h_3.
$$
Here, $b_{ij,k}, i=0,1, j=1,2, k=1,2,\cdots,$ are given by
$$ 
\begin{array}{lll}
&\displaystyle
b_{01,k}=	\frac{ (k\pi)^2+6\cos k\pi +6}{(k\pi)^3},
& 
\displaystyle
b_{11,k}= \frac{2\cos k\pi +4}{(k\pi)^3},\\ 
&\displaystyle
b_{02,k}=\frac{-12 k\pi+12}{(k\pi)^4},
&
\displaystyle
b_{12,k}=\frac{-(k\pi)^2
	-6\cos k\pi +6}{(k\pi)^4} .
\end{array}
$$

{\bf Step 2:}  From proposition \ref{p4.6}, for any fixed $t\in[0,T]$, $v(x,t)\in C^\infty[0,1]$. We have the similar Fourier extension
\begin{equation}\label{vxtexp}
v(x,t)=\alpha_0(t)+\sum_{k=1}^{\infty}\alpha_k(t)\cos k\pi x+\sum_{k=1}^{\infty}\beta_k(t)\sin k\pi x.
\end{equation}
Taking into account that $v(x,0)=0$,  the initial values of are zeros, i.e. 
$$
\alpha_0(0)=0, \;\alpha_k(0)=0, \beta_k(0)=0,\quad \forall k=1,2,\cdots.
$$
Hence, it holds$$
\begin{array}{lcl}
\alpha_0(t)&=&-p_0(t),\\
\alpha_k(t)&=&
\displaystyle
-p_k(t)-i(k\pi)^4\int_0^te^{i(k\pi)^4(t-\tau)}p_k(\tau)d\tau\\
\beta_k(t)&=&
\displaystyle
-q_k(t)-i(k\pi)^4\int_0^te^{i(k\pi)^4(t-\tau)}q_k(\tau)d\tau.
\end{array}
$$
This leads to 
$$
\begin{array}{lll}
v(x,t)&=&
\displaystyle
-\left(p_0(t)+\sum_{k=1}^{\infty}p_k(t)\cos k\pi x+\sum_{k=1}^{\infty}q_k(t)\sin k\pi x\right)\\
&&\displaystyle
-i\sum_{k=1}^{\infty}\int_0^te^{i(k\pi)^4(t-\tau)}
(-12k\pi+12)h_1(\tau)
d\tau\cos k\pi x\\
&&\displaystyle
-i\sum_{k=1}^{\infty}\int_0^te^{i(k\pi)^4(t-\tau)}
(-(k\pi)^2-6\cos k\pi+6) h_3(\tau)
d\tau\cos k\pi x\\
&&\displaystyle 
-i\sum_{k=1}^{\infty}\int_0^te^{i(k\pi)^4(t-\tau)}
((k\pi)^2+6\cos k\pi+6) k\pi h_1(\tau)
d\tau\sin k\pi x\\
&&\displaystyle 
-i\sum_{k=1}^{\infty}\int_0^te^{i(k\pi)^4(t-\tau)}
(2\cos k\pi+4) k\pi h_3(\tau)
d\tau\sin k\pi x.
\end{array}
$$
Taking \eqref{u-v} and \eqref{exfxt1} into account, we arrive at
\eqref{uww} as the forms in \eqref{beta01}-\eqref{betaij}.

Finally,   we let $x'=1-x$,  the situation  can be reduced  to the above cases.  Thus, if $h_1\equiv h_3\equiv 0$ and $h_2(0)=h_4(0)=0, $
$$
u(x,t)=\Big(W_{0,D}h_2+W_{1,D} h_4\Big)\left| _{x\rightarrow1-x}\right. .
$$
The proof is complete.

We now consider the boundary integral induced by the zero-order boundary terms $h_1(t)$ and $h_2(t)$ (recall the definition of $\beta_{ij}$ in \eqref{betaij}):
\begin{eqnarray}
u_{0,h}=W_{0,D} h&\=&
\label{0BC}
\displaystyle S_h+C_h\\
&=&
\displaystyle
\nonumber \sum_{k=1}^{\infty}\int_0^te^{i(k\pi)^4(t-\tau)}
\beta_{01,k}h(\tau)
d\tau\sin k\pi x\\
&&
\displaystyle
\nonumber+
\sum_{k=1}^{\infty}\int_0^te^{i(k\pi)^4(t-\tau)}
\beta_{02,k}h(\tau)
d\tau\cos k\pi x.
\end{eqnarray}


\begin{proposition}\label{0BC-est}
	For any $h\in H^{\frac{3}{4}}(0,T)$, let $u=W_{0,D}(\cdot)h$ and $\Omega_T=(0,1)\times (0,T)$. Then $u$ belongs to $L^4(\Omega_T)\cap C([0,T]; L^2(0,1))$ and satisfies
	$$
	\norm{ u}_{ L^4(\Omega_T)}\leq C\norm{h}_{H^{\frac{3}{4}}(0,T)}
	$$
	and
	$$
	\sup_{0\leq t\leq T}\norm{ u(\cdot, t)}_{L^2(0,1)}\leq C\norm{h}_{H^{\frac{3}{4}}(0,T)}.
	$$
\end{proposition}

{\bf Proof:}  The estimate of $S_h$ is  identical to the one of Proposition \ref{p4.6-1}.  Similarly, one can obtain the estimate of $C_h$. Combining these two estimates we finish the proof.

We now consider the boundary integral induced by the first-order boundary terms $ h_3(t)$ and $h_4(t)$:
\begin{eqnarray}
u_{1,h}=W_{1,D} h
&\=&\label{1BC}
\displaystyle S_h+C_h\\
&=&\nonumber
\displaystyle
-i\sum_{k=1}^{\infty}\int_0^te^{i(k\pi)^4(t-\tau)}
\beta_{11,k} h(\tau)
d\tau\sin k\pi x\\
&&
\nonumber\displaystyle 
-i\sum_{k=1}^{\infty}\int_0^te^{i(k\pi)^4(t-\tau)}
\beta_{12,k} k\pi h(\tau)
d\tau\cos k\pi x.
\end{eqnarray}


\begin{proposition}\label{1BC-est}
	For any $h\in H^{\frac{2}{4}}(0,T)$, let $u=W_{1,D}(\cdot)h$ and $\Omega_T=(0,1)\times (0,T)$. Then $u$ belongs to $L^4(\Omega_T)\cap C([0,T]; L^2(0,1))$ and satisfies
	$$
	\norm{ u}_{ L^4(\Omega_T)}\leq C\norm{h}_{H^{\frac{2}{4}}(0,T)}
	$$
	and
	$$
	\sup_{0\leq t\leq T}\norm{ u(\cdot, t)}_{L^2(0,1)}\leq C\norm{h}_{H^{\frac{2}{4}}(0,T)}.
	$$
\end{proposition}

{\bf Proof:} The estimate of $S_h$ and $C_h$ can be made with the similar techniques as in Proposition \eqref{p4.6-1}. Combining these two estimates we finish the proof. Note that the sharp regularity is given by $C_h$ since $\beta_{12,k}$ has the order $k^2$.

\begin{proposition}\label{p4.6-21}
	Let $0\leq s\leq4$ be given. For any $h\in H_0^{\frac{3+s}{4}}(0,T)$, let $u=W_{0,D} h$. Then $\partial_x^s u$ belongs to $L^4(\Omega_T)\cap C([0,T]; L^2(0,1))$ and satisfies
	$$
	\norm{\partial_x^s u}_{ L^4(\Omega_T)}\leq C\norm{h}_{H^{\frac{3+s}{4}}(0,T)}
	$$
	and
	$$
	\sup_{0\leq t\leq T}\norm{\partial_x^s u(\cdot, t)}_{L^2(0,1)}\leq C\norm{h}_{H^{\frac{3+s}{4}}(0,T)}.
	$$
\end{proposition}

Similar to Proposition \ref{p4.6-21}, it holds
\begin{proposition}\label{p4.6-41}
	Let $0\leq s\leq 4$ be given. For any $h\in H_0^{\frac{2+s}{4}}(0,T)$, let $u=W_{1,D} h$. Then $\partial_x^s u$ belongs to $L^4(\Omega_T)\cap C([0,T]; L^2(0,1))$ and satisfies
	$$
	\norm{\partial_x^s u}_{ L^4(\Omega_T)}\leq C\norm{h}_{H^{\frac{2+s}{4}}(0,T)}
	$$
	and
	$$
	\sup_{0\leq t\leq T}\norm{\partial_x^s u(\cdot, t)}_{L^2(0,1)}\leq C\norm{h}_{H^{\frac{2+s}{4}}(0,T)}.
	$$
\end{proposition}

\section{The nonlinear problem}\label{3.2}

\subsection{With Navier boundary conditions}\label{3.21}

We first set up the well-posed property for the following  problem with Navier boundary conditions:
\begin{equation}\label{4NSDE}
\left\{\begin{array}{ll}
i\partial_t u+  \partial_x^4 u+\lambda |u|^{p-2}u=0, & x\in (0,1),\quad t\in {\mathbb {R}}^+\\
u(x,0)=\varphi(x),&x\in  (0,1)\\
u(0,t)=h_1(t),\quad u(1,t)=h_2(t),\quad &t\in{\mathbb {R}}^+ \\
\partial_x^2 u(0,t)= h_5(t), \quad \partial_x^2 u(1,t)=h_6(t), &t\in{\mathbb {R}}^+
\end{array}\right.
\end{equation}
with $(\varphi, h_1, h_2,  h_5, h_6)\in \mathcal{X}^N_{s,T}$.  Here and thereafter, we denote by
$$
\mathcal{X}^N_{s,T}:=\{(\varphi, h_1, h_2,  h_5, h_6)\Big| \varphi\in H^s(0,1),\;
h_1, h_2\in H_{loc}^{\frac{s+3}{4}}({\mathbb {R}}^+),\;
h_5, h_6\in H_{loc}^{\frac{s+1}{4}}({\mathbb {R}}^+)\}. 
$$

Our first result can be state as:
\begin{theorem}\label{LocalE}
	Let $p\in[3, \infty)$, $s\in (\frac{1}{2}, \frac{9}{2})$ and $\lfloor s\rfloor <p-2$, $T>0$ and $r>0$ be given.  Then there exists a $T^*>0$ such that  the IBVP \eqref{4NSDE} admits a unique solution $u\in C([0,T^*]; H^s(0,1))$, under the following conditions:
	\begin{enumerate}[(1)]
		\item $(\varphi, h_1, h_2,  h_5, h_6)\in \mathcal{X}^N_{s,T}$ with
		$$
		\norm{\varphi}_{H^s(0,1)}+\norm{(h_1,  h_5)}_{X_N}+\norm{(h_2, h_6)}_{X_N} \leq r.
		$$
		\item The  compatibility conditions hold at the corners of the space-time domain, i.e., when $s\in (\frac{1}{2}, \frac{5}{2})$,  $h_1(0)=\varphi(0), h_2(0)=\varphi(1)$; when  $s\in (\frac{5}{2}, \frac{9}{2})$,  $h_1(0)=\varphi(0), h_2(0)=\varphi(1),  h_5(0)=\varphi_{xx}(0), h_6(0)=\varphi_{xx}(1)$.	
	\end{enumerate}
	Moreover, the solution $u$ depends on $(\varphi, h_1, h_2,  h_5, h_6)$ continuously in the corresponding spaces.
\end{theorem}
{\bf Proof:} 
Without loss of generality, we assume that $\varphi(0),\varphi(1),  \varphi_{xx}(0), \varphi_{xx}(1)$ and $ h_i(0), i=1,2,3,4$  equal to $0$.  In fact, the homogenization of  boundary data can be done as follows: Write $ u(x,0)=v(x,0)+\gamma(x)$ with
$$
\gamma(x)=(1-x)(h_1(0)-\frac{1}{6} h_5(0))+\frac{1}{6}(1-x)^3 h_5(0)+x(h_2(0)-\frac{1}{6}h_6(0))+\frac{1}{6}x^3h_6(0).
$$
Then $v$ satisfies homogeneous compatibility conditions and the equation
\begin{equation}\label{modified}
i\partial_t v+\partial_x^4 v+\lambda |v+\gamma|^{p-2}(v+\gamma)=0.
\end{equation}
Similar proof can be given once the estimate of $v$ is done since $\gamma(x)$ is smooth. And we have $$\norm{\gamma}_{H^s(0,1)}\leq C_s,\quad\quad\forall s\geq 0.$$

Note that for $s>\frac{1}{2}$, $H^s(0,1)$ is a Banach algebra, hence, 
\begin{equation}\label{Banach}
\norm{|v+\gamma|^{p-2}(v+\gamma)}_{H^s(0,1)}\leq C\norm{v+\gamma}^{p-1}_{H^s(0,1)}
\end{equation}
with $\lfloor s\rfloor <p-2$.
By Duhamal's principle, we need to solve the integral equation
\begin{eqnarray*}
	v(\cdot, t)&=&W^N(t)\varphi+i\lambda\int_0^tW^N(t-\tau)(|v+\gamma|^{p-2}(v+\gamma))(\cdot,\tau)d\tau\label{umild}\\
	&&+W_{0,N} h_1+(W_{0,N} h_2)|_{x\rightarrow 1-x} +W_{2,N}  h_5+(W_{2,N}  h_5)|_{x\rightarrow 1-x}.
\end{eqnarray*}
Since we are working on the local well-posedness, for fixed $(\varphi, h_1, h_2,  h_5, h_6)\in \mathcal{X}^N_{s,T}$  and $T^*\in(0,T]$,  we set
\begin{eqnarray*}
	\Gamma(v)&=&W^N(t)\varphi+i\chi_{[0,T^*]}(t)\lambda\int_0^tW^N(t-\tau)(|v+\gamma|^{p-2}(v+\gamma))(\cdot,\tau)d\tau \\
	&&+W_{0,N} h_1+(W_{0,N} h_2)|_{x\rightarrow 1-x} +W_{2,N}  h_5+(W_{2,N}  h_5)|_{x\rightarrow 1-x}.
\end{eqnarray*}
For any $v\in C([0,T^*]; H^s(0,1))$, Proposition \ref{p4.6-0}, \ref{p4.6-2}, \ref{p4.6-4} and inequality \eqref{Banach} imply that there exists $C>0$ such that
\begin{equation}\label{Gamestimate}
\norm{\Gamma(v)}_{C([0,\theta]; H^s(0,1))}\leq C(\norm{(\varphi, h_1, h_2,  h_5, h_6)}_{\mathcal{X}_{s,T}^N}+ T^*(\norm{v}_{C([0,\theta]; H^s(0,1))}+C_s)^{p-1}).
\end{equation}
Denote $M=2Cr$. We consider the set $$Y_{s,T^*}:=\{v\in C([0,\theta]; H^s(0,1)),\norm{v}_{C([0,T^*]; H^s(0,1))}\leq M\}.$$
By continuity, we take $T^*$ as small as possible such that $T^* (M+C_s)^{p-1}<\frac{r}{2}$. Here $C$ is the same as in \eqref{Gamestimate}.
Thus,  we have
$$
\norm{\Gamma(v)}_{C([0,\theta]; H^s(0,1))}\leq M
$$
and
$$
\norm{\Gamma(v_1)-\Gamma(v_2)}_{C([0,\theta]; H^s(0,1))}\leq \frac{1}{2}\norm{v_1-v_2}_{C([0,\theta]; H^s(0,1))}
$$
for any $v, v_1, v_2\in Y_{s,T^*}$.
Hence,  $\Gamma$ is a contraction map in $Y_{s,T^*}$. Thus
there exists a unique solution $v\in Y_{s,T^*}$ such that
$$\Gamma(v)=v$$
and $v$ is the unique solution to \eqref{modified} on $[0,T^*]$. We finish the proof of Theorem \ref{LocalE}.
\vskip 5mm

For the well-posedness of the IBVP \eqref{4NSDE} in $H^s(0,1)$ for $s\in [0,\frac{1}{2})$, the space $C([0,T^*]; H^s(0,1))$ is not enough to accomplish the fixed point argument since $H^s(0,1)$ is no longer a Banach algebra. We restrict the argument in $L^4((0,1)\times(0,T^*))\cap C([0,T^*]; H^s(0,1))$. By the same argument above we  have the following lemma:
\begin{lemma}\label{smalls}
	Let $s\in [0, \frac{1}{2})$ and $T>0$. Suppose $r>0$ and $p\in[3,4]$ be given. There exists a $T^*=T^*(r)>0$ such that for any $(\varphi, h_1, h_2,  h_5, h_6)\in \mathcal{X}^N_{s,T}$  with $\norm{(\varphi, h_1, h_2,  h_5, h_6)}_{ \mathcal{X}^N_{s,T}}\leq r$, the IBVP \eqref{4NSDE} admits a unique solution
	$$
	u\in \mathcal{Y}_{s, T^*}:= L^4((0,1)\times(0,T^*))\cap C([0,T^*]; H^s(0,1)),
	$$
	which depends on $(\varphi, h_1, h_2,  h_5, h_6)$ continuously in the corresponding spaces.
\end{lemma}
\textbf{Proof}: First, by repeating the argument in Section 4 of \cite{MR1209299} to set up the local well-posedness in $L^4(\mathbb{T}^2)$. Again, we recall
\begin{eqnarray*}
	\Gamma_\theta(u)&=&W^N(t)\varphi+i\chi_{[0,T^*]}\lambda\int_0^tW^N(t-\tau)(|u|^{p-2}u)(\cdot,\tau)d\tau \\
	&&+W_{0,N} h_1+(W_{0,N} h_2)|_{x\rightarrow 1-x} +W_{2,N}  h_5+(W_{2,N}  h_5)|_{x\rightarrow 1-x}.
\end{eqnarray*}
By Propositions \ref{L4estimates}, \ref{p4.6-1}, \ref{p4.6-3} and noting that $p\in[3,4]$, we have  
$$\norm{\Gamma(u)}_{L^4((0,1)\times(0,T^*))}\leq C\norm{(\varphi,h_1,h_2, h_5,h_6)}_{\mathcal{X}_{0,T}^N}+C{T^*}^{1/5}\|u\|_{L^4((0,1)\times(0,T^*))}^{p-1}$$
and
\begin{eqnarray*}
	&&\hspace{3mm}\norm{\Gamma(u_1)-\Gamma(u_2)}_{L^4((0,1)\times(0,T^*))}\\
	&&\leq C{T^*}^{1/5} (\sum_{i=1,2}\|u_i\|^{p-2}_{L^4((0,1)\times(0,T^*))})
	\norm{u_1-u_2}_{L^4((0,1)\times(0,T^*))}.
\end{eqnarray*}
Let us denote $M=2C\norm{(\varphi,h_1,h_2, h_5,h_6)}_{\mathcal{X}_{0,T}^N}$. If  $0<T^*<T$ is small enough, $\Gamma$ is then a contract map in $\{u\in L^4 \left|\|u\|_{L^4}\leq M\right.\}$. We obtain a unique solution $u\in L^4((0,1)\times(0,T^*))$.  According to Lemma \ref{6.5} and the embedding property of Bourgain space, $u$ is also in $C([0,T^*] H^s(0,1))$ for $s\geq0$.

\subsection{With Dirichlet boundary conditions}

We now consider the problem 
\begin{equation}\label{4-DSDE}
\left\{\begin{array}{ll}
i\partial_t u+  \partial_x^4 u+\lambda |u|^{p-2}u=0, & x\in (0,1),\quad t\in {\mathbb {R}}^+\\
u(x,0)=\varphi(x),&x\in  (0,1)\\
u(0,t)=h_1(t),\quad u(1,t)=h_2(t),\quad &t\in{\mathbb {R}}^+ \\
\partial_x u (0,t)= h_3(t), \quad  \partial_x u (1,t)=h_4(t), &t\in{\mathbb {R}}^+
\end{array}\right.
\end{equation}
with $(\varphi, h_1, h_2,  h_3, h_4)\in \mathcal{X}^D_{s,T}$.  Here and thereafter, we denote by
$$
\mathcal{X}^D_{s,T}:=\{(\varphi, h_1, h_2,  h_3, h_4)\Big| \varphi\in H^s(0,1),\;
h_1, h_2\in H^{\frac{s+3}{4}}([0,T]),\;
h_3, h_4\in H^{\frac{s+2}{4}}([0,T])\}. 
$$
It is  different to the corresponding Navier problem. We first consider the linear initial data problem
\begin{equation}\label{Initial data}
\left\{\begin{array}{ll}
i\partial_t u+  \partial_x^4 u=0, & x\in (0,1),\quad t\in {\mathbb {R}}^+\\
u(x,0)=\varphi(x),&x\in  (0,1),
\end{array}\right.
\end{equation}
From \eqref{LinearIBVPDo} and \eqref{LinearIBVPDe} we get the solution of  \eqref{Initial data} by
$$u(x,t)=u_0(x,t)+u_e(x,t),\quad x\in(0,1), t\in\mathbb R^+.$$
which arises  boundary terms
$$u(0,t)=r_1(t),\quad u(1,t)=r_2(t),\quad  
\partial_x u (0,t)=r_3(t),\quad  \partial_x u (1,t)=r_4(t).
$$
Note that the exact forms of $r_i(t), i=1,2,3,4$ are given by  Lemma \ref{local smoothing}. 

By Proposition \ref{L4estimatesD}, we have
$$
\norm{u}_{C([0,T]; H^s(0,1))}\leq C\norm{\varphi}_{H^s(0,1)}.
$$

Since we are working on local well-posedness, we can restrict $T=1$.
We then  consider the linear problem with zero initial data
\begin{equation}\label{boudary with zero data}
\left\{\begin{array}{ll}
i\partial_t v+  \partial_x^4 v=0, & x\in (0,1),\, t\in [0,1]\\
v(x,0)=0,&x\in  (0,1)\\
v(0,t)=\tilde{h}_1(t),\, v(1,t)=\tilde{h}_2(t),\quad &t\in[0,1] \\
\partial_x v (0,t)= \tilde{h}_3(t), \, \partial_x v(1,t)=\tilde{h}_4(t), &t\in[0,1]
\end{array}\right.
\end{equation}
with $\tilde{h}_i=h_i-r_i, i=1,2,3,4$.
Proposition \ref{BoIn4} shows that
$$v(x,t)=W_{0,D} \tilde{h}_1+W_{1,D}  \tilde{h}_3+\Big(W_{0,D}\tilde{h}_2+W_{1,D} \tilde{h}_4\Big)\left| _{x\rightarrow1-x}.\right.
$$
From propositions \ref{p4.6-21} and \ref{p4.6-41}, for any $s\in [0,4]$, we have
$$\norm{v}_{\mathcal{Y}_{s,T^*}}\leq C(\|(h_1-r_1,h_2-r_2)\|_{H^{\frac{3+s}{4}}(0,1)}+\|(h_3-r_3,h_4-r_4)\|_{H^{\frac{2+s}{4}}(0,1)}).$$
When  $\frac{10}{7}<s\leq \frac{9}{2},$ it can be further bounded by
\begin{equation}\label{boundary condition}
\|h_1\|_{H^{\frac{3+s}{4}}(0,1)}+\|h_2\|_{H^{\frac{3+s}{4}}(0,1)}+\|h_3\|_{H^{\frac{2+s}{4}}(0,1)}+\|h_4\|_{H^{\frac{2+s}{4}}(0,1)}+\|\varphi\|_{H^{s}(0,1)}.\end{equation}
In fact, to get \eqref{boundary condition}, we need to use Lemma \ref{local smoothing} in Appendix. If $0<s\leq1$, we have $\frac{s+3}{4}\leq1$. Thus
$$\|(r_1,r_2)\|_{H^{\frac{s+3}{4}}(0,1)}\leq C\|\varphi_e\|_{H^{\frac{s+3}{8}+\epsilon}(\mathbb T)}.$$ 
If $1<s\leq \frac92$, 
$$\|(r_1,r_2)\|_{H^{\frac{s+3}{4}}(0,1)}\leq C\|\varphi_e\|_{H^{s-\frac12+\epsilon}(\mathbb T)}.$$
If $0<s\leq 2$,
$$\|(r_3,r_4)\|_{H^{\frac{s+2}{4}}(0,1)}\leq C\|\varphi_o\|_{H^{\frac{(s+10)}{8}+\epsilon}(\mathbb T)}.$$ 
If $2\leq s<\frac92$,
$$\|(r_3,r_4)\|_{H^{\frac{s+2}{4}}(0,1)}\leq C\|\varphi_o\|_{H^{s-\frac12+\epsilon}(\mathbb T)}.$$ 
To bound the above terms by $\|\varphi\|_{H^s(0,1)}$, we need 
$$\frac{s+3}{8}<s, \quad s-\frac12<s \quad \text{and} \quad \frac{s+10}{8}<s.$$  
This leads to the lower bound $\frac{10}{7}<s.$

We now consider the following zero initial data and zero  boundary data problem
\begin{equation}\label{final zero initial data and boundary problem}
\left\{\begin{array}{ll}
i\partial_t w+  \partial_x^4 w=|v+u+w|^{p-2}(v+u+w), & x\in (0,1),\quad t\in [0,1]\\
w(x,0)=0,&x\in  (0,1)\\
w(0,t)=0,\quad w(1,t)=0,\quad &t\in [0,1] \\
\partial_x w (0,t)=0, \quad \partial_x w (1,t)=0, &t\in [0,1].
\end{array}\right.
\end{equation}
It is sufficient to prove that \eqref{final zero initial data and boundary problem} has a unique solution in $C([0,T^*]; H^s(0, 1))$ for some $1>T^*>0$. If it is true then $$U(x,t)=u(x,t)+v(x,t)+w(x,t)\in C([0,T^*]; H^s(0, 1))$$ is the unique solution to \eqref{4-DSDE}.
From proposition \ref{p4.6}, we need to solve the following integral equation
$$w=i\int_0^t W^D(t-\tau)|u+v+w|^{p-2}(u+v+w)(\tau)d\tau.$$
We run the fixed point argument again. Let $T^*$   be decided later and $\chi_{[0,T^*]}$ be the characteristic function of the interval $[0,T^*]$.   We denote the map $\Gamma_{\theta} (v)$ by 
$$\Gamma_{T^*} (v)=i\chi_{[0,T^*]}(t)\int_0^t W^D(t-\tau)|u+v+w|^{p-2}(u+v+w)(\tau)d\tau.$$
For $s>10/7$, since  $H^{s}$ is a Banach algebra, we have
$$\|\Gamma_\theta(v)\|_{C([0,T^*]H^s(0,1))}\leq CT^*\|u+v+w\|^{p-1}_{C([0,T^*];H^s(0,1))}$$
with$\lfloor s\rfloor <p-2$. Notice that
$$\|u+v\|_{H^s}\leq C(\|h_1\|_{H^{\frac{3+s}{4}}}+\|h_2\|_{H^{\frac{3+s}{4}}}+\|h_3\|_{H^{\frac{2+s}{4}}}+\|h_4\|_{H^{\frac{2+s}{4}}}+\|\varphi\|_{H^{s}}).
$$
We denote the right side of above inequality by $M$. Then by taking $T^*$ (depending on $M$ and $p$) small enough we have
$$\|\Gamma_{T^*}(v)\|_{C([0,T^*];H^s(0,1))}\leq M.$$
Moreover, by the same argument as in Subsection \ref{3.21}, it holds
$$
\|\Gamma_{T^*}(w_1)-\Gamma_{T^*}(w_2)\|_{C([0,T^*];H^s(0,1))}\leq \frac12\|w_1-w_2\|_{C([0,T^*];H^s(0,1))}.
$$
Hence, we obtain a unique solution to \eqref{4-DSDE} in $C([0,T^*];H^s(0,1))$ for $s>10/7$.

\section {Conclusion}

The study of the IBVP and dynamics of biharmonic nonlinear Schr\"odinger equation with non-homogeneous boundary conditions is important to understand  and design boundary controls when the system evolves in a bounded domain $\Omega$. When dealing with the Navier boundary conditions, the well-posedness can be established via classical PDE techniques in the case of $s>\frac{1}{2}$ in the space $H^s(\Omega)$. In contrast, for $s\leq\frac{1}{2}$, some tools on harmonic analysis are needed and one has to construct some special Banach spaces in order to apply the contraction mapping principle. On the other hand, the well-posedness and mechanism of solutions have been  intensively studied with either the pure initial-value problem posed on the whole domain  $\mathbb{R}^n$   or the periodic boundary initial-value problem posed on the torus $\mathbb{T}^n$.  Boundary integral method provides a systematical way to build a dedicate estimate of the boundary data via similar harmonic analysis techniques. When dealing with the Dirichlet boundary conditions, the situation is much more complicated since  the effects on boundary integrals are not only raised by the boundary data, but also by the initial data. This phenomena restrains $s$ being greater than $10/7$ and it is still an open problem to lower this bound.

As a matter of fact, our results  testify that the stipulation \eqref{KSB} is optimal, which is one of the main contributions of this paper. It is worth mentioning that the optimal regularity of the boundary data is higher than the one expected by the Kato smoothing phenomena and further works need to be done for more general systems.
\section*{Appendix}

\subsection{Fourier series in Dirichlet boundary  problem}

In this section, we first state some facts about Fourier series which we used in dealing with the Dirichlet boundary  problem. Since we concern the mild solutions about the IBVPs, the Fourier series following may not always converge pointwisely but in $L^2(\mathbb T^2)$.

Let $\varphi(x)\in C^4[0,1]$. We make the following two extensions,
$$
\varphi_o(x)=\frac {1}{2}\sgn(x)\varphi(|x|),\quad x\in(-1,1)\quad
\text{with}\quad\varphi_o(0)=\varphi_o(1)=0,$$
and
$$
\qquad  \varphi_e(x)
=\left\{
\begin{array}{ll}
\frac {1}{2} \varphi(|x|), &x\in[-1,0)\cup(0,1] \\
\frac12\varphi(0),&x=0.
\end{array}\right.
$$
Thus $\varphi_0 $ and $\varphi_e$ are $C^4$ functions except the points $\{-1, 0,1\}$. Then we have the following Fourier expansions 
\begin{equation}\label{phioe}
\varphi_o(x)=\sum_{k=1}^{\infty}q_k\sin k\pi x ,
\quad
\varphi_e(x)= p_0
+\sum_{k=1}^{\infty}p_k\cos k\pi x,\quad x\in[-1,1].
\end{equation}
The Fourier coefficients are given by
$$
p_0=\int_0^1\varphi(x)dz,
\quad p_k=\int_0^1\varphi(x)\cos k\pi x dx, 
\quad q_k=\int_0^1\varphi(x)\sin k\pi x dx.
$$
By the Plancherel theorem, we have
$$\|\varphi_o\|_{L^2[-1,1]}^2=\sum_{k=1}^\infty|q_k|^2, \quad\|\varphi_e\|_{L^2[-1,1]}^2=\sum_{k=0}^\infty|p_k|^2.$$
And
$$\|\varphi_o\|_{H^4[-1,1]}^2\sim \sum_{k=1}^\infty|q_k|^2(k\pi)^8, \quad\|\varphi_e\|_{L^2[-1,1]}^2\sim|p_0|^2+\sum_{k=1}^\infty|p_k|^2(k\pi)^8.$$
Actually, for any $0\leq s\leq4$, we have
$$\|\varphi_o\|_{H^s[-1,1]}\sim \sum_{k=1}^{\infty}|q_k|^2(k\pi)^{2s},\quad \|\varphi_e\|_{H^s[-1,1]}^2\sim|p_0|^2+\sum_{k=1}^\infty|p_k|^2(k\pi)^{2s}.$$
From the definitions, it is easy to see that
$$\|\varphi_o\|_{H^s[-1,1]},\|\varphi_e\|_{H^s[-1,1]}\leq C\|\varphi\|_{H^s[0,1]}\leq C\|\varphi\|_{H^s(0,1)},\quad 0\leq s\leq4.$$
We would like to use the notation $\mathbb{T}=\mathbb {R}/2\mathbb {Z}$.

We now consider two functions for $(x,t)\in\mathbb{T}\times\mathbb{R}$:
\begin{equation}\label{uoext}
u_o(x,t)=\sum_{k=1}^{\infty}q_ke^{i(k\pi)^4t}\sin k\pi x ,
\quad
u_e(x,t)= p_0 
+\sum_{k=1}^{\infty}p_k e^{i(k\pi)^4t}\cos k\pi x.
\end{equation} 
Clearly, 
$$
u_o(0,t)=u_o(1,t)=0,  u_o(x,0)=\varphi_o(x).$$ 
For the even part, we have
$$\partial_x u_e(0,t)=\partial_x u_e(1,t)=0,\quad\text{and}\quad u_e(x,o)=\varphi_e(x).$$
Moreover, we have:
\begin{lemma}\label{lema1}
	$u_o(x,t)$ satisfies:
	\begin{equation*}
	\left\{
	\begin{array}{lll}
	&	i\partial_tu_o+\partial_x^4 u_o=0, 
	& (x,t)\in\mathbb{T}\times\mathbb{R},\\
	&%
	u_o(x,0)=\varphi_o(x), 
	& x\in\mathbb{T},\\
	& 
	u_o(0,t)=u_o(1,t)=0,
	&t\in \mathbb{R},\\
	& \displaystyle
	\partial_x u_o(0,t)=\sum_{k=1}^{\infty}q_ke^{i(k\pi)^4t} k\pi\=r_3(t), 
	&
	t\in \mathbb{R},\\
	&\displaystyle
	\partial_x u_o(1,t)=\sum_{k=1}^{\infty}q_ke^{i(k\pi)^4t} k\pi \cos k\pi\=r_4(t), 
	&%
	t\in \mathbb{R}.
	\end{array}
	\right.
	\end{equation*}
	
\end{lemma}
{\bf Proof:} Since $\varphi\in C^4[0,1]$ thus $\varphi_o\in H^4(\mathbb T)$, then we can derivate  the Fourier series in \eqref{uoext} term by term. It is easy to see that $u_0$ satisfies the equation.  $r_3$ and $r_4$ are well defined since the series pointwise converge.
\vskip 5mm

\begin{lemma}\label{lema2} 
	$u_e(x,t)$ satisfies:
	\begin{equation*}
	\left\{
	\begin{array}{lll}
	&
	i\partial_tu_e+\partial_x^4 u_e=0,
	&
	(x,t)\in\mathbb{T}\times\mathbb{R}\\
	& 
	u_e(x,0)=\varphi_e(x),
	&
	x\in \mathbb{T}\\
	&  \displaystyle
	u_e(0,t)=p_0+\sum_{k=1}^{\infty}p_ke^{i(k\pi)^4t}\=r_1(t),	
	&t\in \mathbb{R}\\
	&\displaystyle  u_e(1,t)=p_0+\sum_{k=1}^{\infty}p_ke^{i(k\pi)^4t}\cos k\pi \=r_2(t),	
	&
	t\in \mathbb{R}\\
	&\displaystyle
	\partial_x u_e(0,t)=\partial_x u_e(1,t)=0.
	\end{array}
	\right.
	\end{equation*}
\end{lemma}
{\bf Proof:}  It can be proved in the same way.

So for we introduced the boundary condition $r_1,r_2,r_3$ and $r_4$ for smooth enough function $\varphi\in C^4[0,1]$. We need to extend them into lower regularity spaces. Notice that we satisfies the local mild solutions to the IBVPs, we need only extend $r_1,r_2,r_3,r_4$ into $H^s[0,1]$ for some $0\leq s<4$.  For this aim, we need to set up the following lemma.
\begin{lemma}\label{local smoothing} Let $\varphi$ be as above. $r_1, r_2, r_3, r_4$ are the boundary data raised in Lemma \ref{lema1} and \ref{lema2}. Let $0<s<1$, for any $\epsilon>0$, we have
	\begin{equation}\label{r1r2}
	\|r_1,r_2\|_{H^{s}[0,1]}\leq C \|\varphi_e\|_{H^{\frac12(s+\epsilon)}(\mathbb T)}
	\end{equation}
	and 
	\begin{equation}\label{r3r4}
	\|r_3,r_4\|_{H^{s}[0,1]}\leq C\|\varphi_o\|_{H^{\frac12(s+\epsilon)+1}(\mathbb T)}.
	\end{equation}
	If $1\leq s\leq 2$, we have
	\begin{equation}\label{newr1r2}
	\|r_1,r_2\|_{H^{s}[0,1]}\leq C \|\varphi_e\|_{H^{4(s+\epsilon)-\frac72}(\mathbb T)},
	\end{equation}
	and
	\begin{equation}\label{newr3r4}
	\|r_3,r_4\|_{H^{s}[0,1]}\leq C\|\varphi_o\|_{H^{4(s+\epsilon)-\frac52}(\mathbb T)}.
	\end{equation}
\end{lemma}

{\bf Proof}: Recalling the definition of $r_1(t)$ in Lemma \ref{lema1}, we have
\begin{equation}\label{0-order}
\begin{split}
\|r_1\|_{L^2[0,1]}^2 &= \left|\sum_{k,n\geq0}p_k\bar{p}_n\int_0^1e^{i[(k\pi)^4-(n\pi)^4]t}dt\right| \\
&=\sum_{k}|p_k|^2+\sum_{k\neq n} p_k\bar{p}_n \frac{e^{i[(k\pi)^4-(n\pi)^4]}-1}{i[(k\pi)^4-(n\pi)^4]}\\
&\leq C\|\varphi_e\|^2_{L^2(\mathbb T)}. 
\end{split}
\end{equation}
In the last inequality, we used the Cauchy-Schwarz inequality and use the fact
$$\sum_{k,n\geq 1}\frac{1}{(k+n)^2(k^2+n^2)^2}<\infty.$$
By taking the derivative of $t$ we get
\begin{equation}\label{1-order}
\begin{split}
\|r'_1\|_{L^2[0,1]}^2 &= \left|\sum_{k,n\geq0}i[(k\pi)^4-(n\pi)^4]p_k\bar{p}_n\int_0^1e^{i[(k\pi)^4-(n\pi)^4]t}dt\right| \\
&= \sum_{k\neq n} p_k\bar{p}_n (e^{i[(k\pi)^4-(n\pi)^4]}-1)\\
&\leq C\|\varphi_e\|^2_{H^{\frac12+\epsilon}(\mathbb T)}.
\end{split}
\end{equation}
The last inequality follows from Cauchy-Schwarz and also the fact
$$\sum_{k,n\geq 1}\frac{1}{(kn)^{2(\frac12+\epsilon)}}<\infty.$$
Moreover we also have
\begin{equation}\label{2-order}
\begin{split}
\|r''_1\|_{L^2[0,1]}^2 &= \left|\sum_{k,n\geq 0}-[(k\pi)^8+(n\pi)^8]p_k\bar{p}_n\int_0^1e^{i[(k\pi)^4-(n\pi)^4]t}dt\right| \\
&\leq 2\sum_{k}
(k\pi)^8|p_k|^2+\sum_{k\neq n}[(k\pi)^8+(n\pi)^8] p_k\bar{p}_n \frac{|e^{i[(k\pi)^4-(n\pi)^4]}-1|}{|(k\pi)^4-(n\pi)^4|}\\
&\leq C\|\varphi_e\|^2_{H^{4+\frac12+\epsilon}(\mathbb T)}
\end{split}
\end{equation}
By interpolating between \eqref{0-order} with\eqref{1-order} and \eqref{1-order} with \eqref{2-order}, we obtain
\begin{equation}\label{r1}
\|r_1\|_{H^{s}[0,1]}\leq C \|\varphi_e\|_{H^{\frac12(s+\epsilon)}(\mathbb T)} \;\;\text{for}\,0<s\leq 1
\end{equation}
and
\begin{equation}\label{newr1}
\|r_1\|_{H^{s}[0,1]}\leq C \|\varphi_e\|_{H^{4(s+\epsilon)-\frac72}(\mathbb T)}  \;\; \text{for}\, 1<s\leq2.
\end{equation}
The estimates for $r_2,r_3,r_4$ follow from the same argument. 

Notice that the maps $\varphi_0\to r_3,r_4$ and $\varphi_e\to r_1,r_2$ are all linear maps, we can use the density argument to obtain the boundary conditions $r_1,r_2,r_3$ and $r_4$ for the linear problem with initial data $\varphi\in H^s[0,1]$ for $0\leq s\leq 2$.

\subsection{Some lemmas}

We now introduce  some Lemmas used in the proof of the estimates for the boundary integrals.

\begin{lemma}\label{A-1-1} Let $f\in H^\frac{3}{2}(\mathbb R)$. We have
	$$
	\sum_{k=1}^\infty \|\int_{ 0}^\infty f(\mu) (1-\psi(\mu^4-k^4)) \frac{1}{\mu-k}d\mu\|^2\leq \int_0^\infty (\mu+1)^3|f(\mu)|^2d\mu.
	$$
\end{lemma}
{\bf Proof:} Write
\begin{eqnarray*}
	&&\sum_{k=1}^\infty \|\int_{ 0}^\infty f(\mu) (1-\psi(\mu^4-k^4)) \frac{1}{\mu-k}d\mu\|^2\\
	&&\qquad=\sum_{k=1}^\infty \|\Big (\int_{ 0}^{k-1}+\int_{ k-1}^{k+1}+\int_{k+1}^\infty\Big )f(\mu) (1-\psi(\mu^4-k^4)) \frac{1}{\mu-k}d\mu\|^2\\
	&&\qquad\leq I_1+I_2+I_3.
\end{eqnarray*}
For any fixed $k$ large enough, it is easy to see
$$
\int_0^{k-1}f(\mu)(1-\psi(\mu^4-k^4))\frac{1}{\mu-k}d\mu=\int_0^{k-1}\frac{f(\mu)}{\mu-k}d\mu.
$$
For any $s>\frac12$, by Cauchy-Schwarz, it can be bound by
$$\left(\int_0^{k-1}\frac{d\mu}{(|\mu|+1)^{2s}|\mu-k|^2}\right)^{\frac12}\norm{f}_{H^s(\mathbb R)}.$$
It is element to verify that for $s>\frac12$,
$$\left(\int_0^{k-1}\frac{d\mu}{(|\mu|+1)^{2s}|\mu-k|^2}\right)^\frac12\leq C \frac{1}{k^{\min(s,1)}}.
$$
Summation according to $k$ gives
$$
I_1\leq C \norm{f}_{H^s(\mathbb R)}.
$$
Similarly $I_3$ can also be bounded by $\norm{f}_{H^s(\mathbb R)}$. To study $I_2$, note that in the integrals, the integrand vanishes unless $\mu\geq \sqrt[4]{k^4+1/2}$ or $0\leq \mu\leq \sqrt[4]{k^4-1/2}$. Consequently, for $k$ large enough,
$$\Big (\int_{k-1}^{\sqrt[4]{k^4-1/2}}+\int_{\sqrt[4]{k^4+1/2}}^{k+1}\Big )\frac{d\mu}{|\mu-k|^{2}}\leq C k^3.$$ Then by Cauchy-Schwarz, we have
$$\begin{array}{lll}
&&\displaystyle
\Big (\int_{k-1}^{\sqrt[4]{k^4-1/2}}+\int_{\sqrt[4]{k^4+1/2}}^{k+1}\Big )\|\frac{f(\mu)}{\mu-k}\|d\mu \\
&\leq &\displaystyle
C k^{\frac{3}{2}}\Big (\int_{k-1}^{k+1}|f(\mu)|^2d\mu\Big )^{\frac12}\approx \Big (\int_{k-1}^{k+1}|\mu|^3|f(\mu)|^2d\mu\Big )^{\frac12}.
\end{array}$$
Then we finish the proof by summation up $k$.
\vskip 5mm
\begin{lemma}\label{p2.1}
	The set $\{(k, k^4), \,\| \; n\in\mathbb{Z}\}$ has bounded $\Lambda$-constant, i.e.
	$$
	\norm{\sum_{k\in \mathbb{Z}}a_k e^{i(kx+k^4t)}}_{L^4({\mathbb {T}}^2)}\leq c\Big (\sum_{k+\mathbb{Z}}|a_k|^2\Big )^{1/2}
	$$
\end{lemma}
{\bf Proof:} We write
$$
\norm{\sum_{k\in \mathbb{Z}}a_k e^{i(kx+k^4t)}}_{L^4({\mathbb {T}}^2)}=\norm{\sum_{k,l\in\mathbb{Z}}a_k\bar{a}_l e^{i[(k-l)x+(k^4-l^4)t]}}_{L^2({\mathbb {T}}^2)}.
$$
By Plancherel theorem, it equals to
$$\Bigg\{\sum_{\xi,\eta\in\mathbb{Z}}\Big(\sum_{(k,l)\in A(\xi,\eta)}a_k\bar{a}_l\Big)^2\Bigg\}^{\frac{1}{2}}.$$
Here $A(\xi,\eta)=\{(k,l)\in\mathbb{Z}^2;k-l=\xi,k^4-l^4=\eta\}.$ It is not hard to proof that $A(\xi,\eta)\cap A(\xi^\prime,\eta^\prime)=\emptyset.$ And meanwhile, for any fixed $\xi,\eta\in\mathbb{Z}$,
$$\#A(\xi,\eta)\leq 3.$$
We finish the proof by Cauchy-Schwarz.
\vskip 5mm
\begin{lemma}\label{L4E}
	We have  that
	$$
	\norm{f}^2_{L^4({\mathbb {T}}^2)}\leq C\Big [\sum_{m,n\in\mathbb{Z}}\Big (|n-m^4|+1\Big )^{2s}|\hat f(m, n)|^2\Big ]^{1/2}
	$$
	holds for any any $s>\frac{1}{4}$.
\end{lemma}
{\bf Proof:} The main idea of the proof follows \cite{MR1209299}. We here give a proof for the convenience of the readers. We write
$$\norm{f}^2_{L^4({\mathbb {T}}^2)}=\norm{f\bar{f}}_{L^2({\mathbb {T}}^2)}.$$
By Plancherel theorem, it equals to
$$\sum_{m,n\in \mathbb{Z}^2}\Big|\sum_{m_1,n_1\in \mathbb{Z}^2}\hat{f}(m_1,n_1)\hat{\bar{f}}(m-m_1,n-n_1)
\Big|^2.$$
For  $s>\frac14$, we insert the $(|n-m^4|+1)^s$ and use Cauchy-Schwarz to bound it by
$$
\sup_{m,n\in \mathbb{Z}^2}
A(m,n)\sum_{m,n\in\mathbb{Z}}\Big (|n-m^4|+1\Big )^{2s}|\hat f(m, n)|^2
$$
with
$$
A(m,n )=\sum_{m_1,n_1\in \mathbb{Z}^2}\frac{1}{(|n_1-m_1^4|+1)^{2s}(|n-n_1-(m-m_1)^4|+1)^{2s}}.
$$
We suffer to show
\begin{equation}\label{6.2}
\sup_{m,n\in \mathbb{Z}^2}A(m, n)\leq C
\end{equation}
with $s>\frac14$. Now we fix $m, n, m_1$ and denote $\tilde{n}=n_1-m_1^4$, the summation in \eqref{6.2}
can be written as
$$
\sum_{m_1,{\tilde{n}}\in \mathbb{Z}^2}\frac{1}{(|\tilde{n}|+1)^{2s}(|\tilde{n}-a|+1)^{2s}}
$$
with
$$a=a(n,m,m_1)=n-m_1^4-(m-m_1)^4.$$
For any $j\geq0$, we say that $m_1\in A_j$ if for fixed $n,m\in\mathbb{Z}$ s.t.
$$m_1\in \mathbb{Z},\quad 2^{j-1}\leq |a|<2^{j},\text{for} j\geq1,\quad \text{or}\quad |a|<1,\quad\text{for}\quad j=0.$$
Then \eqref{6.2} can be bounded by
$$
\sum_{j\geq0\atop m_1\in A_j}\sum_{\tilde{n}\in\mathbb{Z}}\frac{1}{(|\tilde{n}|+1)^{2s}(|\tilde{n}-a|+1)^{2s}}\leq C\sum_{j\geq0}
\frac{\#(A_j)}{2^{2js}}.
$$
We claim that \begin{equation}\label{6.4}\#(A_j)\leq C 2^{\frac{j}{2}}.\end{equation} With this claim and note that $s>\frac14$, we can finish the proof. We suffer to set up \eqref{6.4}. First if
$|(m-m_1)^3-m_1^3|\geq 2^{j/2}.$ We use to mean value theorem to obtain that
$$\#(A_j)\leq 2^{j/2}.$$
On the other hand, the mean value theorem again shows that there are at most $2^{j/2}$ many $m_1\in\mathbb{Z}$ such that $|(m-m_1)^3-m_1^3|<2^{j/2}.$
We finish the proof.
\vskip 5mm
As a direct consequence of Lemma \ref{L4E}, it holds
\begin{lemma}\label{6.5}
	We have
	$$
	\norm{f}_{L^4({\mathbb {T}}\times{\mathbb {R}})}\leq c\Big [\int_{{\mathbb {R}}}\sum_{m\in\mathbb{Z}}\Big (|\lambda-m^4|+1\Big )^{2s}|\hat f(m, \lambda)|^2d\lambda\Big ]^{1/2}
	$$
	with $s>\frac14$.
\end{lemma}
{\bf Proof:} It can be done by the same argument of Lemma \ref{L4E}. 
\vskip 5mm
\begin{lemma}\label{p37-middle}
	We have  that
	$$
	\|\sum_{k=\lfloor\sqrt[4]{2\lambda}\rfloor}^\infty \sin(k\pi x)\frac{1}{\sqrt[4]k-{\lambda}}\|\leq C|x|^{\alpha-1}\frac{1}{(1+\sqrt[4]{\lambda})^{1-\alpha}},
	$$
	with $\alpha\in(\frac{3}{4}, 1)$.
\end{lemma}
{\bf Proof:} Let
$$
S_n=\sum_{k=1}^n\sin k\pi x=\frac{\sin((n+1)\pi x/2)\sin(n\pi x/2)}{\sin(\pi x/2)}, \quad \hbox{for}\; n=1,2,\cdots.
$$
For any $\alpha\in[0,1]$ and $0<x\leq 1$, $|S_n| \leq C \frac{n^\alpha}{|x|^{1-\alpha}}$
\footnote{It is due to the fact that  $|S_n| \leq C \frac{|nx|^2}{|x|}\leq C \frac{|nx|^\alpha}{|x|}$ when $(n+1)x\leq 1$, and  $|S_n| \leq C \frac{1}{|x|}\leq C \frac{|nx|^\alpha}{|x|}$  when $(n+1)x> 1$.}
.
Consequently, for any $n\geq 2$,
\begin{eqnarray*}
	\sum_{k=\lfloor\sqrt[4]{2\lambda}\rfloor}^n\frac{1}{k-\sqrt[4]{\lambda}}(S_k-S_{k-1})
	=
	\sum_{k=\lfloor\sqrt[4]{2\lambda}\rfloor}^n\frac{1}{k-\sqrt[4]{\lambda}}S_k-\sum_{k=\lfloor\sqrt[4]{2\lambda}\rfloor}^n\frac{1}{k-\sqrt[4]{\lambda}}S_{k-1}\\
	=\sum_{k=\lfloor\sqrt[4]{2\lambda}\rfloor}^{n-1}\Big (\frac{1}{k-\sqrt[4]{\lambda}}-\frac{1}{k+1-\sqrt[4]{\lambda}}\Big )S_k+\frac{1}{n-\sqrt[4]{\lambda}}S_n-\frac{1}{\lfloor\sqrt[4]{2\lambda}\rfloor -\sqrt[4]{\lambda}}S_{\lfloor\sqrt[4]{2\lambda}\rfloor-1}.
\end{eqnarray*}
Choose $3/4<\alpha<1$ and let $n\rightarrow \infty$, we have
$$
\begin{array}{lll}
\displaystyle\|\sum_{k=\lfloor\sqrt[4]{2\lambda}\rfloor}^\infty\frac{1}{k-\sqrt[4]{\lambda}}\sin k\pi x\|
&\leq &
\displaystyle C|x|^{\alpha-1}\Big (\frac{\lambda^{\alpha/4}}{\sqrt[4]{\lambda}+1}+\sum_{k=\lfloor\sqrt[4]{2\lambda}\rfloor}^\infty\frac{k^\alpha}{(k-\sqrt[4]{\lambda})^2}\Big )\\
&\leq &
\displaystyle C|x|^{\alpha-1}\Big (\frac{\lambda^{\alpha/4}}{\sqrt[4]{\lambda}+1}+\sum_{k=1}^\infty\frac{1}{(k+\sqrt[4]{\lambda})^{2-\alpha}}\Big )\\
&\leq &
\displaystyle  \frac{C|x|^{\alpha-1}}{(\sqrt[4]{\lambda}+1)^{1-\alpha}}.
\end{array}
$$
Here is the conclusion of the Fourier series of $\sin(\sqrt i a(\pi-x))$.
\vskip 5mm
\begin{lemma}\label{expansion-s}
	we have  that
	$$
	\sum_{k=1}^{\infty}\frac{k^3+ika^2}{k^4+a^4}\sin (kx)=\frac{\pi}{2}\frac{\sin(\sqrt i a(\pi-x))}{\sin(\sqrt i a\pi)}.
	$$
\end{lemma}

{\bf Proof:} Denote by $\varphi_k=\int_0^\pi \sin(\sqrt i a(\pi-x))\sin (kx) dx$. Then
$$
\begin{array}{lll}
\varphi_k
& = &
\displaystyle -\frac{1}{k}\cos (kx)\sin(\sqrt i a(\pi-x))\Big|_0^\pi\\
&    &
\displaystyle
+\frac{1}{k}\int_0^\pi (-\sqrt i a)\cos(\sqrt i a(\pi-x))\cos (kx) dx
\\
& = &
\displaystyle \frac{1}{k}\sin(\sqrt i a\pi)+\frac{ia^2}{k^2}\int_0^\pi \sin(\sqrt i a(\pi-x))\sin (kx) dx
\\
& = &
\displaystyle  \frac{1}{k}\sin(\sqrt i a\pi)+\frac{ia^2}{k^2}\Big (\frac{1}{k}\sin(\sqrt i a\pi)+\frac{ia^2}{k^2}\varphi_k\Big ),
\end{array}
$$
which means for any $k=1,2,\cdots, $
$$
\varphi_k=\frac{k^3+ika^2}{k^4+a^4}\sin(\sqrt i a\pi).
$$
The result is straightforward since for any $x\in(0,\pi)$, the Fourier series of $\sin(\sqrt i a(\pi-x))$ is given by
$$
\sin(\sqrt i a(\pi-x))=\sum_{k=1}^{\infty}\frac{2}{\pi}\varphi_k\sin (kx).
$$
\begin{lemma}\label{sina}
	For $\sin(\sqrt i a)$, it holds
	$$
	\sin(\sqrt i a)=\frac{e^{-\frac{a}{\sqrt 2}}e^{\frac{a}{\sqrt 2}i}-e^{\frac{a}{\sqrt 2}}e^{-\frac{a}{\sqrt 2}i}}{2i}
	$$
\end{lemma}
The optimality of the parameters are given in the following results.
\begin{lemma}\label{opt0}
	All regularities of the boundary terms are optimal. More precisely, we have
	\begin{itemize}
		\item $h\in H^{\frac{3}{4}}(0,T)$  is optimal for the zero-order boundary data $h_1$ and $h_2$;
		\item $h\in H^{\frac{2}{4}}(0,T)$  is optimal for the first-order boundary data $h_3$ and $h_4$;
		\item $h\in H^{\frac{1}{4}}(0,T)$  is optimal for the second-order boundary data $h_5$ and $h_6$.
	\end{itemize}
\end{lemma}
{\bf Proof:} We give the proof for $h_1$ and $h_2$.  The rest can be verified with the same methodology.

Note that Proposition \ref{p4.6-1}  implies that
$$
\norm{u_{0,h}}_{L^2((0,1)\times(0,T))}\leq C_T\norm{u_{0,h}}_{H^{\frac{3}{4}}(0,T))},
$$
with the definition
$$
u_{0,h}= \sum_{k=-\infty}^\infty \beta_k e^{i(k\pi)^4t}\int_0^t e^{-i(k\pi)^4\tau}h(\tau)d\tau  e^{i k\pi x}.
$$
Assume that $h(t)$ has the Fourier series
$$
h(t)=\sum_{n=-\infty}^{\infty} e^{i n\pi^4 t} a_n \qquad \hbox{with}\qquad a_n=\int_{0}^{\frac{2\pi}{\pi^4}}e^{-i n\pi^4 t} h(t)dt.
$$
It follows that
$$
\begin{array}{lll}
u_{0,h}&=&\displaystyle
\sum_{k=-\infty}^\infty -(k\pi)^3e^{i(k\pi)^4t+ik\pi x}\sum_{n=-\infty}^{\infty} \int_0^t e^{i(n-k^4)\pi^4\tau}a_n d\tau \\
&=&\displaystyle
\sum_{k=-\infty}^\infty -(k\pi)^3e^{i(k\pi)^4t+ik\pi x}
\Big (\sum_{n\neq k^4} \frac{ e^{i(n-k^4)\pi^4t-1}}{i(n-k^4)\pi^4}a_n+t a_{ k^4}\Big ) \\
&=&\displaystyle
\sum_{k=-\infty}^\infty -(k\pi)^3e^{ik\pi x}
\sum_{n\neq k^4} \frac{ e^{in\pi^4t}- e^{ik^4\pi^4t}}{i(n-k^4)\pi^4}a_n\\
&  &\displaystyle
+\sum_{k=-\infty}^\infty -(k\pi)^3e^{i(k\pi)^4t+ik\pi x}t a_{ k^4}.
\end{array}
$$
Choose $h(t)$ such that
$$
a_{ k^4}=\int_{0}^{\frac{2\pi}{\pi^4}}e^{-i k\pi^4 t} h(t)dt=0, \qquad  n\in \mathbb{Z}.
$$
The last term of the above formula vanishes and we have
$$
\begin{array}{lll}
u_{0,h} 
&  =  &
\displaystyle
\sum_{k=-\infty}^\infty -(k\pi)^3e^{ik\pi x}
\sum_{n\neq k^4} \frac{ e^{in\pi^4t}a_n}{i(n-k^4)\pi^4}\\
&      &
\displaystyle
+
\sum_{k=-\infty}^\infty (k\pi)^3e^{ik\pi x+ik^4\pi^4t}
\sum_{n\neq k^4} \frac{  a_n}{i(n-k^4)\pi^4}.
\end{array}
$$
Since $e^{ik\pi x+in\pi^4t}$ and $e^{ik\pi x+ik^4\pi^4t}$ are orthogonal as $n\neq k^4$, it holds that
\begin{equation}\label{contra3}
\begin{array}{lll}
&&\displaystyle
\norm{u_{0,h}}^2_{L^2((0,1)\times(0, \frac{2}{\pi^3}))}\\
&=&\displaystyle
\sum_{k=-\infty}^\infty
\Big (
\sum_{n\neq k^4}(k\pi)^6  \frac{  a_n^2}{(n-k^4)^2\pi^8}
+
\Big (\sum_{n\neq k^4}(k\pi)^3  \frac{  a_n}{(n-k^4)\pi^2}\Big )^2
\Big ) \\
&\geq&\displaystyle
\pi^{-2}  \sum_{k=-\infty}^\infty k^6a_{k^4+1}^2.
\end{array}
\end{equation}
We now prove the optimality of the regularity by contradiction. Suppose that there exists a constant $C$ such that
\begin{equation}\label{contra1}
\norm{u_{0,h}}^2_{L^2((0,1)\times(0, \frac{2}{\pi^3}))}\leq C \norm{h}^2_{H^{\alpha}(0,\frac{2}{\pi^3})}
\end{equation}
for some $\alpha\in(0,\frac{3}{4})$.  We choose $h_n$  as the form
\begin{equation}\label{contra2}
h_n(t)=\sum_{k\neq0}^{|k|\leq n}\frac{1}{|k|^\beta}e^{i(k^4+1)\pi^4t}, \qquad \hbox{for }\; n=1,2,\cdots
\end{equation}
with $\beta\in(\frac{1+8\alpha}{2},\frac{7}{2})$. Indeed, it is due to the fact that  the inequality
$$
\sum_{k\neq0}\|\frac{(k^4+1)^\alpha}{|k|^\beta}\|^2<\infty
$$
holds whence $2\beta-8\alpha>1$. Consequently $h_n$ belongs to $H^{\alpha}(0,\frac{2}{\pi^3})$ as $n$ tends to infinity.

However, by taking the boundary data $h_n$ as in \eqref{contra2}, combining \eqref{contra3} and \eqref{contra1},  we arrive at
$$
C \norm{h_n}^2_{H^{\alpha}(0,\frac{2}{\pi^3})}
\geq \norm{u_{0,h_n}}^2_{L^2((0,1)\times(0, \frac{2}{\pi^3}))}
\geq \pi^{-2}  \sum_{k=-\infty}^\infty k^6a_{k^4+1}^2
= \sum_{k\neq 0}^{|k|\leq n} \frac{1}{\pi^2|k|^{2\beta-6} }.
$$ 
The last term of the above formula tends to infinity as $n$ tends to infinity. This is a contradiction and the proof is complete.

\section*{References}

\bibliography{mybibfile}

\end{document}